\documentclass[11pt,english]{article}
\usepackage{geometry}
\usepackage{authblk}
\usepackage{setspace}
\usepackage{graphicx}
\usepackage{amsmath,color}
\usepackage{mathbbol}
\usepackage{multirow}
\newcommand{\tabincell}[2]{\begin{tabular}{@{}#1@{}}#2\end{tabular}}
\usepackage{cases}
\usepackage{subfig}

\newtheorem{remark}{Remark}

\title{A deterministic solver for the linear Boltzmann model of a single mono-directional proton beam}

\author[1$\dag$]{Xiaojiang Zhang}
\author[2$\ddag$]{Xuemin Bai}
\author[3$\dag$,*]{Min Tang}

\affil[1]{School of Mathematics and Statistics, Central South University, Changsha 410083, China}
\affil[2]{ Mevion Medical Equipment Co., Ltd.,
  135 Yuanfeng Road, Kunshan City, Jiangsu Province, China}
  \affil[3]{School of Mathematical Sciences, Institute of Natural Sciences, MOE-LSC, Shanghai Jiao Tong University, Shanghai, 200240, P.R. China.}
\affil[*]{Corresponding author. Email: tangmin@sjtu.edu.cn}

\date{}

\begin{document}
	
	\maketitle
	
	%%%%%% Abstract %%%%%%
	\begin{abstract}
		The linear Boltzmann model for proton beams is a six-dimensional partial differential equation (PDE). We propose a deterministic solver for the linear Boltzmann model based on scattering decomposition and depth-splitting methods. The main idea is to first divide the protons into primary protons and scattering protons, whose equations are derived using the source iteration method. We then treat depth as the time variable in classical time-evolutionary problems and apply the depth-splitting method. In the depth-splitting method, the full operator is decomposed into three parts, with each subsystem being easily parallelizable, which is crucial for efficient simulations. The resulting discretization exhibits second-order convergence in both the depth and energy variables. The dose distributions obtained from our solver are compared with those from Monte Carlo simulations for various materials and heterogeneous cases.
	\end{abstract}
	
	%%%%%% Main Text %%%%%%
	
	\section{Introduction}
	 Radiation therapy can be divided into different types according to the different external beams of particles \cite{uilkema2012proton}. Photon therapy is the most mature and widely used radiation treatment, while proton therapy attracts more and more attention since it has much less side effects.
	%Radiation is typically administered using external beam photon therapy. On the other hand, proton therapy, in existence for over 60 years but initially limited to research laboratories until the 1990s, has experienced rapid clinical growth, with over 50 facilities worldwide \cite{levin2005proton}.
	The physical mechanisms of photon therapy and proton therapy are significantly different \cite{paganetti2013biological}.  In the photon therapy, the dose distribution as a function of depth displays   a maximum dose close to the entrance, then followed by an exponentially decreasing energy deposition. On
	the other hand, due to the mass of a proton  being 1832
	times greater than an electron,  protons mainly undergo small-angle scattering events. So protons   deposit maximum energy  close to the end of their range in the proton therapy, resulting in a 'Bragg' peak. Moreover, the proton dose is small at  the entrance and exit of the beam, and these properties make protons especially suited for therapy.
	
	Simulating the dose distributions is crucial in  radiation treatment. There are two traditional methods 
	for calculating the dose profile of a proton beam:  pencil-beam algorithms \cite{hong1996pencil,pedroni2005experimental} and the Monte Carlo method \cite{perl2012topas,paganetti2008clinical}. The pencil-beam method is based on the the Fermi-Eyges theory \cite{eyges1948multiple}, and this method is computationally efficient, but the  accuracy is low, especially when heterogeneity of the human body is considered \cite{szymanowski2002two}. On the other hand, the Monte Carlo (MC) method simulates individual interacting protons and can achieve a
	satisfactory accuracy. However, despite continuous research efforts to  accelerate MC methods, their  computational costs are high and  suffer from unavoidable statistical fluctuations.
	
	In order to achieve a computationally feasible model with comparable accuracy, the linear Boltzmann model is considered as a useful alternative to describe the motion of the protons \cite{gifford2006comparison,bedford2019calculation,uilkema2012proton}. Compared to the pencil-beam method, the full linear Boltzmann model provides better accuracy and is applicable to both homogeneous and heterogeneous materials. On the other hand, the accuracy is comparable to the MC method, while there exists no statistical noise. The main difficulty of simulating the linear Boltzmann model for radiation therapy is the high dimensionality. The system involves six independent variables, including the position in  physical space, energy, and proton movement direction in phase space, which results in a high computational cost. Since the required computational time is crucial in real clinical applications, to design an efficient PDE solver, it is important to reduce the computational cost and make it easy to parallelize.

	Linear Boltzmann solver has been successfully developed and used for photon therapy. The Acuros algorithm of Vassiliev et al was adopted and further developed by Varian for the Eclipse treatment planning system \cite{vassiliev2010validation}.
	Linear Boltzmann solver for proton therapy is less mature and no well established fast solver is available yet.
	Simplified model for proton transport is considered in \cite{burlacu2023deterministic} for the laterally homogeneous, in-depth heterogeneous geometry and only the primary protons are considered. In order to reduce the computational cost, the authors construct an approximation by separation of variables. However, depth-dependent mean energy has to be introduced to determine the model parameter, which introduces additional errors \cite{geback2012analytical}. 

 %The Fermi pencil-beam model is fast but allows only for the case $\sigma_{tr}=\sigma_{tr}(x)$, i.e., for layered heterogeneities. Its application  to nonlayered media is difficult. However, we use numerical scheme to treat $\sigma_{tr}$, which is not involved in Fermi–Eyges  theory. Thus our scheme can handle the case $\sigma_{tr}=\sigma_{tr}(x,y,z)$.
	
	% Thanks to  a lot of deterministic solvers that have been developed for neutron transport \cite{larsen2006overview,lawrence1986progress,ciolini2002simplified,talamo2013numerical},  an efficient and accurate numerical method to the linear Boltzmann equation can be achieved. In terms of angular discretization, deterministic methods can be roughly divided into two types:  the discrete-ordinates ($S_N$) method \cite{bgc1} and the spherical harmonics ($P_N$) method \cite{Frank-Klar-Larsen-2007}. The discrete ordinate method discretizes the  distribution function along several specific angular directions and  the  energy density is reconstructed using the  quadrature rule.  In the $S_N$ method, the solution remains positive, however, it can exhibit non-physical artifacts, referred to as ray-effects. In the framework of the $P_N$ method, the distribution function in the angular space is approximated using a series expansion of polynomials. While the $P_N$ method does not maintain the positivity of the solution and can potentially result in oscillatory approximations,  efficient numerical treatment of scattering terms can be obtained. 
	
	In this paper, we propose a deterministic solver based on scattering decomposition and depth splitting  method for the linear Boltzmann model which describes the distribution of single mono-directional proton beam. The main idea is to firstly divide the protons into primary protons and scattering protons. In order to solve the governing equation of each proton type, the depth is then treated as the time in classical time-evolutionary problems. Strang depth splitting method is then employed to solve the depth evolutionary problem, in which the full operator is decomposed into three parts: the transport operator in space, the scattering operator in velocity and the continuous slowing down operator in energy. There are two main advantages of the proposed method:  (1) When different values of the variables of other dimensions are used, the subsystem can be solved in parallel. Therefore, each subsystem is easy to parallel which allows for a fast simulation. This is important for real clinical applications; (2) Each subsystem is solved by the second order Crank-Nicolson (CN) method. Due to the good property of the CN method, the implicit part can be solved explicitly and the computational cost can be reduced by allowing for a large grid size in depth. 

According to the numerical results, the accuracy of our proposed method in terms of dose distribution is comparable to the classical Monte Carlo method and meets clinical standards. Moreover, the numerical solutions exhibit no noise, making them more accurate in spatial regions where there are fewer protons.
On the other hand, compared to the semi-analytical method proposed in \cite{burlacu2023deterministic}, our method can handle laterally heterogeneous cases and accounts for catastrophic scatter interactions with a nucleus, thus being more suitable for clinical applications.
In terms of computational efficiency, the proposed deterministic solver can be effectively parallelized. However, we have only conducted a preliminary test of the scheme's performance in this work; a full efficiency comparison with Monte Carlo is beyond the scope of the current paper.
	
	The organization of this paper is as follows. Section 2 introduce the linear Boltzmann model for the proton, and simplified approximations to the collision integrals. 
	The detailed energy discretization, the depth splitting method
	and
	fully discretized system are presented in Sections 3. 
	In Section 4, we compare the obtained dose distribution  with a general MC code, FLUKA. Different materials and heterogeneous cases are considered and some indexes that are interesting for medical physical community are presented. A very good agreement of the linear Boltzmann solver with the MC results is achieved for different cases. Finally the paper is concluded with some discussions in Section 5. 
	
	%\textcolor{red}{Why pronton therapy becomes popular, its advantage, simulation based on MC,pencil beam, disadvantage, talk about the main stream software, advantage of deterministic method(no noise, because pronton almost move straight forward, the material geometry is less important,  neutron transport developed a lot of deterministic solver, uncertainty analysis); difficulty of the deterministic method(high dimensionality, review the modeling papers, their methods have to solve a large system). our contribution, propose the splitting method for primary pronton, Advantage: 1)different material(homogeneous and heterogeneous), different energy(works for injected energy has a distribution), 2)carefully design the splitting to make the decoupling of different dimensions, the system is  easy to parallel,  }

	% %%%%%% Figures %%%%%%

	\section{The transport model for the protons}
	
	There are three main types of interactions for protons in matter: inelastic Coulomb scatter with atomic electrons, elastic Coulomb scatter with atomic nuclei, and catastrophic scatter interactions with nuclei \cite{newhauser2015physics,enferadi2017nuclear}. During these interactions, the energy of protons is deposited, and the movement direction of protons changes, and secondary particles can be generated. Most protons travel nearly straight in matter due to the mass of the proton being much larger than the electron. During its frequent inelastic interactions with atomic electrons, the kinetic energy of the protons is lost continuously. When a proton passes close to a nucleus whose mass is larger than the proton, it experiences elastic Coulomb scattering and deflects from its original trajectory. In catastrophic scatter interactions, the proton enters the nucleus, and the nucleus will emit secondary particles. The catastrophic scatter interactions between nuclei and protons are less frequent and remove the incident proton from the beam. Thus, the primary protons are those that only involve inelastic Coulomb scatter with atomic electrons and elastic Coulomb scatter with atomic nuclei. Secondary particles originate from catastrophic scatter interactions \cite{paganetti2002nuclear,zhang2022primary}, whose scattering kernel is more complex, and we focus on reconstructing the scattering kernel from MC code in this paper.

	The linear Boltzmann equation for the probability density distribution of the protons writes \cite{asadzadeh2018priori, borgers1996asymptotic,burlacu2023deterministic}:
	\begin{equation}\label{fp}
		\begin{aligned}
			\boldsymbol{\Omega}\cdot\nabla\psi(\textbf{r},\boldsymbol{\Omega},E)&=\int_0^\infty\sigma_{i,s}(\textbf{r},\boldsymbol{\Omega},E+W\to E)\psi(\textbf{r},\boldsymbol{\Omega},E+W)\text{d}W-\sigma_{i,t}(\textbf{r},E)\psi(\textbf{r},\boldsymbol{\Omega},E)\\
		&+\int_{4\pi}\sigma_{e,s}(\textbf{r},\boldsymbol{\Omega'}\to\boldsymbol{\Omega},E)\psi(\textbf{r},\boldsymbol{\Omega'},E)\text{d}\boldsymbol{\Omega'}-\sigma_{e,t}(\textbf{r},E)\psi(\textbf{r},\boldsymbol{\Omega},E)\\
  &+\int_{4\pi}\int_E^\infty\sigma_{c,s}(\textbf{r},\boldsymbol{\Omega'}\to\boldsymbol{\Omega},E'\to E)\psi(\textbf{r},\boldsymbol{\Omega'},E')\text{d}\boldsymbol{\Omega'}\text{d}E'-\sigma_{c,t}(\textbf{r},E)\psi(\textbf{r},\boldsymbol{\Omega},E), 
		\end{aligned}
	\end{equation}
	where $\psi(\textbf{r},\boldsymbol{\Omega},E)$ is the probability density distribution of protons moving in direction $\boldsymbol{\Omega}$ at position $\textbf{r}$  with energy $E$ ; $\textbf{r}=(x,y,z)$ is the space vector;  $\boldsymbol{\Omega}=(\mu,\eta,\xi)$ is the unit vector on a unit sphere; $\sigma_{i,s}$, $\sigma_{e,s}$ and $\sigma_{c,s}$ denote the scattering cross section for inelastic scattering, elastic scattering and catastrophic scattering respectively; $\sigma_{i,t}$, $\sigma_{e,t}$ and $\sigma_{c,t}$ represent the total  cross section for inelastic scattering, elastic scattering and catastrophic scattering respectively, which can be obtained by taking the integral of the scattering  cross section with respect to $\boldsymbol{\Omega}$; $W$ is defined as the energy loss during the interaction process.

	As discussed in \cite{zheng1993overview}, the angular deflection during the inelastic interaction with atomic electrons and energy loss during the Coulomb elastic scatter interactions are negligible.
	The first two terms on the right-hand side of equation \eqref{fp} represents the protons with energy higher than $E$ that lose some energy and become protons with energy $E$ after the inelastic interaction with atomic electrons, while their movement directions do not change.   The third and fourth terms on the right-hand side of equation \eqref{fp} represent the protons with movement direction $\boldsymbol{\Omega}'$ change their movement direction to $\boldsymbol{\Omega}$ after the elastic interaction with nuclei, while their energy remains unchanged. Moreover, the last two terms stands for all the particles with energy $E$ and moving in direction $\Omega$ following a catastrophic scatter interaction with a nucleus.

Equation \eqref{fp} provides the probability density distribution of the  protons, whose integral can be used to calculate the dose distribution due to the protons. However, $\psi(\textbf{r},\boldsymbol{\Omega},E)$ is a $6$-dimensional integral differential equation, standard PDE solvers are too expensive and can be very slow. For clinical applications, the solver has to be efficient and easy to parallelize. The nonlocal integral terms are not favorable for efficiency, thus, in \cite{burlacu2023deterministic,pomraning1992fokker,zheng1993overview}, the following approximations are proposed for the collision integrals: 
	
	%\subsection{ Boltzmann-Fokker-Planck approximation}
	\begin{itemize}
		\item The continuous slowing-down approximation (CSDA) and the energy-loss straggling (ELS) approximation are applied to the inelastic scattering. Since the energy loss during a single inelastic interaction between the proton and the atomic electron is very small, a continuous energy loss with stopping power  and straggling coefficient  can be used to approximate the  inelastic stopping process effectively. More precisely,
		\begin{equation}
			\begin{aligned}
				& \int_0^\infty\sigma_{i,s}(\textbf{r},\boldsymbol{\Omega},E+W\to E)\psi(\textbf{r},\boldsymbol{\Omega},E+W)\text{d}W-\sigma_{i,s}(\textbf{r},E)\psi(\textbf{r},\boldsymbol{\Omega},E)\\
				\approx&\rho(\textbf{r})\frac{\partial}{\partial E}\left(S(\textbf{r},E)\psi\right)+\rho(\textbf{r})\frac{1}{2}\frac{\partial^2}{\partial E^2}\left(T(\textbf{r},E)\psi\right),
			\end{aligned}
		\end{equation}
		where $S(\textbf{r},E)$ and $T(\textbf{r},E)$ are stopping power and straggling coefficient respectively; $\rho(\textbf{r})$ is the target weight density. The coefficients $S(\textbf{r},E)$ and $T(\textbf{r},E)$ are respectively the first and second moments of cross section $\sigma_{i,s}$, which can viewed as the mean and the deviation around the mean during the continuous energy loss process. Its derivation can be found in \cite{uilkema2012proton}, for the convenience of readers, we have put them in the Appendix~\ref{apen01} and \ref{apen02}.

		\item The Fokker-Planck approximation for the elastic scatter process.  Let 
		\begin{equation}\label{defi_}
			\mu=\text{cos}\vartheta,\quad\eta=\sqrt{1-\mu^2}\text{cos}\phi,\quad \xi=\sqrt{1-\mu^2}\text{sin}\phi,
			\quad (\vartheta,\phi)=[0,\pi]\times[0,2\pi],\end{equation}
		with $\vartheta$ and $\phi$ being respectively the polar and azimuth angle. For forward peaked small-angle scattering, the Fokker-Planck approximation has been derived by the asymptotic analysis in \cite{pomraning1992fokker}. In \cite{hensel2006deterministic}, this approximation has been applied to the elastic scatter process between proton and atomic electrons. More precisely, 
		\begin{equation}\label{app_fp}
			\begin{aligned}
				&   \int_{4\pi}\sigma_{e,s}(\textbf{r},\boldsymbol{\Omega'}\to\boldsymbol{\Omega},E)\psi(\textbf{r},\boldsymbol{\Omega'},E)\text{d}\boldsymbol{\Omega'}-\sigma_{e,s}(\textbf{r},E)\psi(\textbf{r},\boldsymbol{\Omega},E)\\
				\approx &\rho\frac{\sigma_{tr}}{2}\left[\frac{\partial}{\partial\mu}(1-\mu^2)\frac{\partial}{\partial\mu}+\frac{1}{1-\mu^2}\frac{\partial^2}{\partial\phi^2}\right]\psi,
			\end{aligned}
		\end{equation}
		where $\sigma_{tr}$ is the momentum transfer cross section, which is defined by
		\begin{equation}\label{eq:mt}
			\sigma_{tr}(E)=2\pi\int_{-1}^1(1-\mu)\sigma_{e,s}(E,\mu)d\mu.
		\end{equation}
		The momentum transfer cross section $\sigma_{tr}(E)$ represents the average angular deviation of the protons per unit travelled distance, whose specific
		expression is given in Appendix A.3.   
	\end{itemize}
	
Then the linear Boltzmann transport equation \eqref{fp}  can be rewritten as:
	\begin{equation}\label{f-p1}
 \begin{aligned}
     	\mu\frac{\partial\psi}{\partial x}+\eta\frac{\partial\psi}{\partial y}+\xi\frac{\partial\psi}{\partial z}&=\rho\frac{\sigma_{tr}}{2}\left[\frac{\partial}{\partial\mu}(1-\mu^2)\frac{\partial}{\partial\mu}+\frac{1}{1-\mu^2}\frac{\partial^2}{\partial\phi^2}\right]\psi+\rho\frac{\partial}{\partial E}\left(S(\textbf{r},E)\psi\right)\\ &+\rho\frac{1}{2}\frac{\partial^2}{\partial E^2}\left(T(\textbf{r},E)\psi\right)+\int_{4\pi}\int_E^\infty\sigma_{c,s}\psi\text{d}\boldsymbol{\Omega'}\text{d}E'-\sigma_{c,t}(\textbf{r},E)\psi.
      \end{aligned}
	\end{equation}
Let
	$$
	u=\frac{\eta}{\mu},\quad v=\frac{\xi}{\mu},
	$$

 In actual clinical applications, 
for mono-directional proton beam from the nozzle, the directional and energy distribution of proton emission are highly concentrated. On the proton radiotherapy machine, the energy and angle distribution of the emitted protons follow a Gaussian distribution, with a standard deviation of 0.35\ MeV for energy and $0.13^{\circ}$ for angle. Moreover, owing to the mass of a proton being much larger than the mass of an electron (1832 times the mass), most protons in proton therapy mainly undergo small-angle scattering events and travel in a straight line in the depth. Therefore,
  we consider only an approximation that is associated with Fermi-Eyges theory. This approximation assumes that the transverse cosines are small, thus it is inaccurate once the beam has sufficiently expanded in the transverse direction. We will show later from numerical tests that the energy deposition of this approximation can achieve the accuracy requirement for clinical applications.

In the Fermi-Eyges approximation, the moving directions of all protons in the injected proton beam are centered in direction $(\mu,\eta,\xi)=(1,0,0)$.
 The proton beam is sharply peaked near $\mu=1$ under the small angle approximation. 
	Noting the definitions in \eqref{defi_}, from $\mu\approx 1$, 
	$$
	\eta\ll 1,\quad \xi\ll 1.
	$$
	% As in \cite{borgers1996asymptotic}, let 
	% $$
	% u=\frac{\eta}{\mu},\quad v=\frac{\xi}{\mu}.
	% $$
The Fokker-Planck operator in brackets on the right hand side of equation \eqref{f-p1} can be replaced by the Laplacian operator:
	\begin{equation}\label{eq:L}
	L=\frac{\partial^2}{\partial u^2}+\frac{\partial^2}{\partial v^2}.
	\end{equation}
	Then one can consider $\psi(x,y,z,\mu,u,v,E)=\hat\psi (x,y,z,1,u,v,E)$, and the linear Boltzmann transport equation \eqref{f-p1} can be rewritten into an equation for $\hat{\psi}$. For simplicity, we drop the hat and consider the following equation:
		\begin{equation}\label{f-e}
 \begin{aligned}
     	\frac{\partial\psi}{\partial x}+u\frac{\partial\psi}{\partial y}+v\frac{\partial\psi}{\partial z}&=\rho\frac{\sigma_{tr}}{2}\left(\frac{\partial^2}{\partial u^2}+\frac{\partial^2}{\partial v^2}\right)\psi+\rho\frac{\partial}{\partial E}\left(S(\textbf{r},E)\psi\right)+\rho\frac{1}{2}\frac{\partial^2}{\partial E^2}\left(T(\textbf{r},E)\psi\right)\\ &+\int_{4\pi}\int_E^\infty\sigma_{c,s}\psi\text{d}u'\text{d}v'\text{d}E'-\sigma_{c,t}(\textbf{r},E)\psi.
      \end{aligned}
	\end{equation}
%In the machine of Mevion, the standard deviation of the directional distribution is .... Moreover, in our numerical examples, the numerical results of the model \eqref{f-e} match those of Monte Carlo very well. Therefore, in order to reduce the amount of calculation, we will consider the discretization of model \eqref{f-e} next.
 %However, the authors never discuss the implications of this approximation on the accuracy of the solution. In effect, the authors use an approximation that can be described as Fermi-Eyges theory combined with continuous-slowing down and Fokker-Planck energy straggling. This is not as accurate as the full Fokker-Planck approximation. It is well-known that Fermi-Eyges theory gives inaccurate results after significant expansion of the beam. Why do the authors' computational results not reflect this?

In  \eqref{f-e}, the injected proton in proton therapy enter at position $x=0$, and mainly move in the direction perpendicular to the $YZ$ plane. The  phase-space density of injected protons can have a distribution in $E$,$u$,$v$,$y$, and $z$, which indicates that the injected protons allow different energies, positions and movement directions, as far as the included angles between the movement directions and the direction perpendicular to the $YZ$ plane are small. Moreover, this approximation allows one to deal with a simplified two-dimensional Laplacian instead of the full Fokker-Planck operator as in \eqref{app_fp}. The diffusion coefficient $\sigma_{tr}$ is independent of $(u,v)$, thus one can either use semi-analytical solutions as in the Fermi-Eyges theory or classical fast solvers for two-dimensional Laplacians.

 In this paper, we focus on the design of a deterministic algorithm for this single mono-directional proton beam.
The beams considered in the numerical tests are of a diameter appropriate for proton therapy.
As will be observed in the numerical tests, the dose obtained by simulating \eqref{f-e} is accurate enough when compared to the results obtained by MC methods. The $\rho$ in \eqref{f-e} can depend on $(x,y,z)$, thus the model in \eqref{f-e} can simulate cases where both lateral and in-depth distributions are heterogeneous.
 Given that the specific expression of the catastrophic scattering cross sections are not publicly available, we use data from FLUKA to fit them, which can be found in the Appendix~\ref{sec-d}. By using the fitted catastrophic scattering cross section, we can meet the clinical accuracy requirements in numerical experiments. 

	\section{The numerical scheme}
	In this section, we will introduce an efficient numerical discretization for the transport model \eqref{f-e}.  Since in the most clinically relevant cases, the injected protons have almost the same energy level, which indicates that, the energy distribution at $x=0$ is very concentrated. Thus the energy discretization should allow discontinuities in the solution,  we adopt the DG  method for energy discretization. Moreover, due to the specific set up in proton therapy, we consider only one particular type of boundary conditions, the protons are only injected at $x=0$ and zero elsewhere at the boundary of the spatial domain. Due to this specific set up, one can use splitting method along the $x$-axis which is efficient and easy to parallel. Finally, since the included angles between the proton movement directions and the direction perpendicular to the $YZ$ plane are small, we consider $\psi(x,y,z,u,v,E)=0$ for $(u,v)$ at the boundary of the $UV$ domain for $\forall x,y,z,E$. The sizes of the $UV$ domain depend on the material. As we can observe in the numerical examples, for typical material considered in proton therapy, water, bone and air, it is enough to take $(u,v)\in[-1,1]\times[-1,1]$. Due to the zero boundary condition and the Laplacian operator in the velocity variable $u,v$,  we use Fourier spectral method in the velocity variable $(u,v)$. In the subsequent part, we will give the details of the depth splitting method.
 %%%%%%%%%%%%%%%%%%%%%%%%%%%%%%%%%%%%%%%%
\subsection{Scattering decomposition method}

The proton in the therapy can be further divided into the primary protons and scattering protons \cite{hklt}. The primary protons are interacted with the material by the ionization process, and the scattering protons are generated by inelastic scattering. Therefore, we divide the distribution of the protons into the primary  part  and scattering part such that
\begin{equation}
\psi=\psi_p+\sum_{k=1}^\infty\psi_{s,k}.\label{scatterdecompose}
\end{equation}
Here $\psi_p$ is the distribution of the primary proton, and $\psi_{s,1}$ is the distribution of protons that  have catastrophic scattered once, and $\psi_{s,k}$ ($k=2,3,\cdots$) is the distribution of  protons that  have  catastrophic scattered $k$times . Moreover, $\psi_p$, $\psi_{s,1}$ and $\psi_{s,k}$ ($k=2,3,\cdots$) satisfy the following equations:
\begin{subequations}\label{eqn:pri-se}
\begin{numcases}{}
\frac{\partial\psi_p}{\partial x}+u\frac{\partial\psi_p}{\partial y}+v\frac{\partial\psi_p}{\partial z}=\rho\frac{\sigma_{tr}}{2}\left(\frac{\partial^2}{\partial u^2}+\frac{\partial^2}{\partial v^2}\right)\psi_p+\rho\frac{\partial}{\partial E}\left(S(\textbf{r},E)\psi_p\right)-\sigma_{c,t}(\textbf{r},E)\psi_p \nonumber\\
   \hspace{3.8cm}+\rho\frac{1}{2}\frac{\partial^2}{\partial E^2}\left(T(\textbf{r},E)\psi_p\right), \label{eqn:pri-se1}\\
			\frac{\partial\psi_{s,1}}{\partial x}+u\frac{\partial\psi_{s,1}}{\partial y}+v\frac{\partial\psi_{s,1}}{\partial z}=\rho\frac{\sigma_{tr}}{2}\left(\frac{\partial^2}{\partial u^2}+\frac{\partial^2}{\partial v^2}\right)\psi_{s,1}+\rho\frac{\partial}{\partial E}\left(S(\textbf{r},E)\psi_{s,1}\right)-\sigma_{c,t}(\textbf{r},E)\psi_{s,1}\nonumber\\\hspace{3.8cm}+\rho\frac{1}{2}\frac{\partial^2}{\partial E^2}\left(T(\textbf{r},E)\psi_{s,1}\right) +\int_{4\pi}\int_E^\infty\sigma_{c,s}\psi_p\text{d}u'\text{d}v'\text{d}E',\label{eqn:pri-se2}\\
   	\frac{\partial\psi_{s,k}}{\partial x}+u\frac{\partial\psi_{s,k}}{\partial y}+v\frac{\partial\psi_{s,k}}{\partial z}=\rho\frac{\sigma_{tr}}{2}\left(\frac{\partial^2}{\partial u^2}+\frac{\partial^2}{\partial v^2}\right)\psi_{s,k}+\rho\frac{\partial}{\partial E}\left(S(\textbf{r},E)\psi_{s,k}\right)-\sigma_{c,t}(\textbf{r},E)\psi_{s,k}\nonumber\\\hspace{3.8cm}+\rho\frac{1}{2}\frac{\partial^2}{\partial E^2}\left(T(\textbf{r},E)\psi_{s,k}\right) +\int_{4\pi}\int_E^\infty\sigma_{c,s}\psi_{s,k-1}\text{d}u'\text{d}v'\text{d}E'.\label{eqn:pri-se3}
\end{numcases}
\end{subequations}
It is easy to check that, if the expansion on the right hand side of \eqref{scatterdecompose} converges, $\psi_p+\sum_{k=1}^\infty\psi_{s,k}$ provide the solution to the original equation \eqref{f-e}.
This decomposition is analog to the source iteration method \cite{mathews2009adaptive,kusch2023robust} and we have checked numerically in subsection 4.1 that for the proton therapy, $\psi=\psi_p+\psi_{s,1}$ provides an energy depletion that is accurate enough. \eqref{eqn:pri-se1} describes the dynamics of primary protons that can be solved independently. Then, primary protons undergo catastrophic scattering provide the source term (integral term) for protons which have scattered once whose dynamics are described by \eqref{eqn:pri-se2}.

Except the source terms, the operator is the same for $\psi_p$ and $\psi_k$ ($k=1,2,\cdots$). As far as one can get an efficient solver for $\psi_p$, one can use the same solver to get $\psi_k$ ($k=1,2,\cdots$) by taking into account the source terms.  In the subsequent part we focus on equation  \eqref{eqn:pri-se1}, and for ease of writing, use $\psi$ instead of $\psi_p$.
 
	\subsection{Energy discretization} 
	For the energy discretization, the multi-group method is the most popular one \cite{mihalas1999foundations,pomraning2005equations}. In the multi-group method, the continuous energy space  is divided into $G$ groups, and each energy interval is denoted by $(E_{g-1/2},E_{g+1/2})$ ($g=1,\cdots, G$) with $g=G$ being the highest energy group. In the discontinuous Galerkin (DG) method, the fluxes inside each energy groups are assumed to be continuous, and the solution at the group faces is allowed to be discontinuous. Thus, similar to \cite{burlacu2023deterministic,uilkema2012proton},
	we adopt DG method to allow for the discontinuous distribution in energy and to guarantee energy conservation. 
	
	In each energy group, the probability density distribution function $\psi$ is approximated by a polynomial of order $N_p$:
	\begin{equation}\label{eqn:df_p}
		\psi(x,y,z,u,v,E)=\sum_{l=1}^{N_p}\psi_g^l(x,y,z,u,v)\phi_g^l(E),
	\end{equation}
	where $\phi_g^l(E)$ are the basis functions in the interval $K=(E_{g-1/2},E_{g+1/2})$. Let $\Delta E_g=E_{g+1/2}-E_{g-1/2}$, the basis functions satisfy the orthogonally property:
	\begin{equation}\label{eq:orthogonal}
		\int_{E_{g-1/2}}^{E_{g+1/2}}\phi_g^i\phi_g^j\text{d}E=\frac{\Delta E_g}{2j-1}\delta_{i,j}.
	\end{equation} 
	By substituting \eqref{eqn:df_p} into \eqref{f-e}, one gets
	\begin{equation}\label{eq:group}
		\begin{aligned}
			&\sum_{l=1}^{N_p}\phi_g^l\partial_x\psi_g^l+\sum_{l=1}^{N_p}\phi_g^l u\partial_y\psi_g^l+\sum_{l=1}^{N_p}\phi_g^l v\partial_z\psi_g^l+\sum_{l=1}^{N_p}\phi_g^l\sigma_{c,t}\psi_g^l\\
			=
			&\sum_{l=1}^{N_p}\phi_g^l\rho\frac{\sigma_{tr}}{2}\left(\frac{\partial^2}{\partial u^2}+\frac{\partial^2}{\partial v^2}\right)\psi_g^l+\sum_{l=1}^{N_p}\rho\partial_E\left(S(\textbf{r},E)\psi_g^l\phi_g^l\right)+\sum_{l=1}^{N_p}\rho\frac{1}{2}\frac{\partial^2}{\partial E^2}\left(T(\textbf{r},E)\psi_g^l\phi_g^l\right).
		\end{aligned}
	\end{equation}
	We use a second order DG in this paper. Let $N_p=2$,  and the basis functions are
	$$
	\phi_g^1(E) = 1, \qquad \phi_g^2(E) = \frac{2}{\Delta E_g}(E-E_g),
	$$
	with $E_g=\frac{E_{g+1/2}+E_{g-1/2}}{2}$.
	As discussed in \cite{uilkema2012proton}, the stopping power $S(\mathbf{r},E)$ plays a crucial role in the transport equation as it  serves as the primary factor influencing energy deposition. Moreover, the value of $S(\mathbf{r},E)$ changes fast with respect to $E$ (see Fig.~\ref{fig:s_p}), thus we approximate the stopping power within each energy group by a linear function such that 
	\begin{equation}\label{eq:appSE}
		\begin{aligned}
			S(E)&=\frac{S_{g+\frac{1}{2}}+S_{g-\frac{1}{2}}}{2}+ \frac{S_{g+\frac{1}{2}}-S_{g-\frac{1}{2}}}{\Delta E_g}(E-E_g),\\
			&=\frac{S_{g+\frac{1}{2}}+S_{g-\frac{1}{2}}}{2} \phi^1_g(E)+\frac{S_{g+\frac{1}{2}}-S_{g-\frac{1}{2}}}{2} \phi^2_g(E),
		\end{aligned}
	\end{equation}
	where $S_{g \pm \frac{1}{2}}=S(E_{g \pm \frac{1}{2}})$, and we define $T(\mathbf{r},E)$ as the same way. In the subsequent part, we use the orthogonality in \eqref{eq:orthogonal} and test \eqref{eq:group} by $\phi_g^1$ and $\phi_g^2$ to get the evolutionary equations for $\psi_g^1$ and $\psi_g^2$:

	\begin{itemize}
		\item Testing the term $\sum_{l=1}^2\partial_E\big(S(\mathbf{r},E)\psi_g^l\phi_g^l\big)$ by $\frac{\phi_g^1}{\Delta E_g}$ yields:
		$$
		\begin{aligned}
			&\frac{1}{\Delta E_g} \int_{E_{g-\frac{1}{2}}}^{E_{g+\frac{1}{2}}} \frac{\partial}{\partial E} \big(S(E) \psi\big) \phi_g^1 \text{d}E\\
			=&\frac{1}{\Delta E_g}\left[(S(E) \psi) \phi_g^1\right]\Big|^{E_{g+\frac{1}{2}}}_{E_{g-\frac{1}{2}}}-\frac{1}{\Delta E_g}\int_{E_{g-\frac{1}{2}}}^{E_{g+\frac{1}{2}}}  S(E) \psi\frac{\partial}{\partial E} \phi_g^1 \text{d}E\\
			\approx &\frac{S_{g+\frac{1}{2}}}{\Delta E_g}\left(\psi_{ g+1}^1-\psi_{ g+1}^2\right)-\frac{S_{g-\frac{1}{2}}}{\Delta E_g}\left(\psi_{ g}^1-\psi_{ g}^2\right),
		\end{aligned}
		$$
		where, in the second equality, $(S(E) \psi) \phi_g^1\Big|_{E_{g+\frac{1}{2}}}=S(E_{g+\frac{1}{2}}) \psi(E_{g+\frac{1}{2}})$. According to \eqref{eqn:df_p}, $\psi(E_{g+\frac{1}{2}})$ has two different values at the group interface: $\psi_g^1+\psi_g^2$ and $\psi_{g+1}^1-\psi_{g+1}^2$, the upwind value is used to approximate $\psi(E_{g+\frac{1}{2}})$.
		\item By using \eqref{eq:appSE}, one can test the term $\sum_{l=1}^2\partial_E\big(S(\mathbf{r},E)\psi_g^l\phi_g^l\big)$ by $\frac{3\phi_g^2}{\Delta E_g}$ and get:
		$$
		\begin{aligned}
			&\frac{3}{\Delta E_g} \int_{E_{g-\frac{1}{2}}}^{E_{g+\frac{1}{2}}} \frac{\partial}{\partial E} (S(E) \psi) \phi_g^2 \text{d}E\\
			=&\frac{3}{\Delta E_g}\left[(S(E) \psi) \phi_g^2\right]\Big|^{E_{g+\frac{1}{2}}}_{E_{g-\frac{1}{2}}}-\frac{3}{\Delta E_g}\int_{E_{g-\frac{1}{2}}}^{E_{g+\frac{1}{2}}}  S(E) \psi\frac{\partial}{\partial E} \phi_g^2 \text{d}E\\
			=&\frac{3}{\Delta E_g}\left[(S(E) \psi) \phi_g^2\right]\Big|^{E_{g+\frac{1}{2}}}_{E_{g-\frac{1}{2}}}-\frac{6}{(\Delta E_g)^2}\int_{E_{g-\frac{1}{2}}}^{E_{g+\frac{1}{2}}}  \left(\frac{S_{g+\frac{1}{2}}+S_{g-\frac{1}{2}}}{2} \phi^1_g(E)+\frac{S_{g+\frac{1}{2}}-S_{g-\frac{1}{2}}}{2} \phi^2_g(E)\right)\psi\text{d}E\\
			\approx &3 \frac{S_{g+\frac{1}{2}}}{\Delta E_g}\left(\psi^1_{g+1}-\psi^2_{g+1}\right)-3 \frac{S_{g-\frac{1}{2}}}{\Delta E_g}\left(\psi^1_{g}-\psi^2_{g}\right)-3 \frac{S_{g+\frac{1}{2}}+S_{g-\frac{1}{2}}}{\Delta E_g} \psi^1_{g}-\frac{S_{g+\frac{1}{2}}-S_{g-\frac{1}{2}}}{\Delta E_g} \psi^2_{g}. \\
		\end{aligned}
		$$
		Here  in the last approximation, we have used \eqref{eqn:df_p} and the orthogonally property of the basis functions.
		\item For the energy straggling
		operator $\frac{\partial^2}{\partial E^2}\left(T(\textbf{r},E)\psi\right)$, we firstly divide it into two parts: $\frac{\partial}{\partial E}\left(T\left(\frac{\partial \psi}{\partial E}\right)\right)$ and $\frac{\partial}{\partial E}\left(\left(\frac{\partial T}{\partial E}\right)\psi\right)$, and the second part can be treated like the stopping power operator, moreover,   one can use the symmetric interior penalty Galerkin method \cite{müller2017symmetric} to treat the first part, and define the jump and the average at the edges of an energy interval:
  $$
  \left[\psi\right]=\psi^--\psi^+,\quad \{\psi\}=\frac{\psi^-+\psi^+}{2},
  $$
  where $\psi^-=\psi(E^-)$ and $\psi^+=\psi(E^+)$.
  One can test the first part  by $\frac{\phi_g^1}{\Delta E_g}$ and get
		$$\begin{aligned}
			&\frac{1}{\Delta E_g} \int_{E_{g-\frac{1}{2}}}^{E_{g+\frac{1}{2}}} \phi_g^1  \frac{\partial}{\partial E}\left(T\frac{\partial \psi}{\partial E}\right)\text{d}E\\
   =&\frac{1}{\Delta E_g}\int_{\partial K}\left(T \partial_E\psi \right)n_K\cdot\phi_g^1\text{d}s-\frac{1}{\Delta E_g}\int_{E_{g-\frac{1}{2}}}^{E_{g+\frac{1}{2}}} T\frac{\partial \psi}{\partial E}\frac{\partial \phi_g^1}{\partial E}  \text{d}E-\frac{\alpha}{\Delta E_g}\left[\psi\right]\left[\phi_g^1\right]\\  
   =&\frac{1}{\Delta E_g}\left(\{T\partial_E\psi\}\cdot[\phi_g^1]+[\psi]\cdot\{T\partial_E\phi_g^1\}\right)-\frac{1}{\Delta E_g}\int_{E_{g-\frac{1}{2}}}^{E_{g+\frac{1}{2}}} T\frac{\partial \psi}{\partial E}\frac{\partial \phi_g^1}{\partial E}  \text{d}E-\frac{\alpha}{\Delta E_g}\left[\psi\right]\left[\phi_g^1\right]
   :=Q_1,
  % \approx &\frac{2}{\Delta E_g}\frac{T_{g+1} +T_{g} }{2} \frac{ \psi^1_{g+1}- \psi^1_g}{\Delta E_g+\Delta E_{g+1}}-\frac{2}{\Delta E_g}\frac{T_{g-1} +T_{g} }{2} \frac{\psi^1_g-\psi^1_{g-1}}{\Delta E_g+\Delta E_{g-1}},
		\end{aligned}$$
		where 
  %$T_{g \pm 1}=T(E_{g \pm 1})$, and
  the term $\frac{\alpha}{\Delta E_g}\left[\psi\right]\left[\phi_g^1\right]$ is the penalty term, where the penalty parameter $\alpha$ is defined as
  $$
  \alpha=\frac{[max(deg(\phi_g^l))+1]^2}{2},
  $$
  with $deg(\phi_g^l)$ is the polynomial degree  of the basis functions.
%\textcolor{red}{3. The authors claim to use a discontinuous Galerkin approximation in energy, which is fine for the continuous-slowing-down operator, but it is problematic for the straggling operator. A straightforward FEM approach is not possible because the derivative of the distribution function at the highest group energy is a delta-function and thus unbounded. This unbounded quantity appears in the weak form. For instance, it appears in the equation above line 182, but the authors nonetheless appear to arrive at a simple cell-centered discretization for the straggling term in that equation. I don't see where this discretization can possibly be consistent with a discontinuous Galerkin approximation.}
  
		\item Finally, by testing the term  $\frac{\partial}{\partial E}\left(T\left(\frac{\partial \psi}{\partial E}\right)\right)$ by $\frac{3\phi_g^2}{\Delta E_g}$, one can get:
		$$
		\begin{aligned}
			&\frac{3}{\Delta E_g} \int_{E_{g-\frac{1}{2}}}^{E_{g+\frac{1}{2}}} \phi_g^2 \frac{\partial}{\partial E}\left(T\frac{\partial \psi}{\partial E}\right)\text{d}E
			\\=&\frac{3}{\Delta E_g}\int_{\partial K}\left(T \partial_E\psi \right)n_K\cdot\phi_g^2\text{d}s-\frac{3}{\Delta E_g}\int_{E_{g-\frac{1}{2}}}^{E_{g+\frac{1}{2}}} T\frac{\partial \psi}{\partial E}\frac{\partial \phi_g^2}{\partial E}  \text{d}E-\frac{3\alpha}{\Delta E_g}\left[\psi\right]\left[\phi_g^2\right]\\
			 =&\frac{3}{\Delta E_g}\left(\{T\partial_E\psi\}\cdot[\phi_g^2]+[\psi]\cdot\{T\partial_E\phi_g^2\}\right)-\frac{3}{\Delta E_g}\int_{E_{g-\frac{1}{2}}}^{E_{g+\frac{1}{2}}} T\frac{\partial \psi}{\partial E}\frac{\partial \phi_g^2}{\partial E}  \text{d}E-\frac{3\alpha}{\Delta E_g}\left[\psi\right]\left[\phi_g^2\right]:=Q_2.
		%	\approx &\frac{6}{\Delta E_g} \frac{T_{g+1} \psi^1_{g+1}-T_g \psi^1_g}{\Delta E_g+\Delta E_{g+1}}+\frac{6}{\Delta E_g} \frac{T_g \psi^1_g-T_{g-1} \psi^1_{g-1}}{\Delta E_g+\Delta E_{g-1}}-\frac{6T_{g+\frac{1}{2}}}{(\Delta E_g)^2}\left(\psi_{ g+1}^1-\psi_{ g+1}^2\right)+\frac{6T_{g-\frac{1}{2}}}{(\Delta E_g)^2}\left(\psi_{ g}^1-\psi_{ g}^2\right), \\
		\end{aligned}
		$$
	%	with $T_{g\pm\frac{1}{2}}=\frac{T_{g}+T_{g\pm 1}}{2}$.
	\end{itemize}
	In summary, one can get the following system for $\psi_g^1$ and $\psi_g^2$:
	\begin{equation}\label{eqn:0041-d}
		\left\{
		\begin{aligned}
			&\frac{\partial\psi_g^1}{\partial x}+u\frac{\partial\psi_g^1}{\partial y}+v\frac{\partial\psi_g^1}{\partial z}=\rho\frac{\sigma_{tr,g}}{2}\left(\frac{\partial^2}{\partial u^2}+\frac{\partial^2}{\partial v^2}\right)\psi_g^1+\rho\frac{S_{g+\frac{1}{2}}}{\Delta E_g}\left(\psi_{ g+1}^1-\psi_{ g+1}^2\right)-\rho\frac{S_{g-\frac{1}{2}}}{\Delta E_g}\left(\psi_{ g}^1-\psi_{ g}^2\right)\\
			&\qquad-\sigma_{c,t,g}\psi^1_g+\frac{\rho}{2}\frac{T_{g+1}-T_g}{(\Delta E_g)^2}\left(\psi_{ g+1}^1-\psi_{ g+1}^2\right)-\frac{\rho}{2}\frac{T_{g}-T_{g-1}}{(\Delta E_g)^2}\left(\psi_{ g}^1-\psi_{ g}^2\right)+\frac{\rho}{2} Q_1, \\
			&\frac{\partial\psi_g^2}{\partial x}+u\frac{\partial\psi_g^2}{\partial y}+v\frac{\partial\psi_g^2}{\partial z}=\rho\frac{\sigma_{tr,g}}{2}\left(\frac{\partial^2}{\partial u^2}+\frac{\partial^2}{\partial v^2}\right)\psi_g^2+3\rho \frac{S_{g+\frac{1}{2}}}{\Delta E_g}\left(\psi^1_{g+1}-\psi^2_{g+1}\right)
			-3 \rho\frac{S_{g-\frac{1}{2}}}{\Delta E_g}\left(\psi^1_{g}-\psi^2_{g}\right)\\
			&\qquad
			-3\rho \frac{S_{g+\frac{1}{2}}+S_{g-\frac{1}{2}}}{\Delta E_g} \psi^1_{g}-\rho\frac{S_{g+\frac{1}{2}}-S_{g-\frac{1}{2}}}{\Delta E_g} \psi^2_{g}-\sigma_{c,t,g}\psi^2_g+\frac{3\rho}{2} \frac{T_{g+1}-T_{g}}{(\Delta E_g)^2}\left(\psi^1_{g+1}-\psi^2_{g+1}\right)\\
			&\qquad-\frac{3\rho}{2} \frac{T_{g}-T_{g-1}}{(\Delta E_g)^2}\left(\psi^1_{g}-\psi^2_{g}\right)-\frac{3\rho}{2} \frac{T_{g+1}-T_{g-1}}{(\Delta E_g)^2}\psi^1_{g}+\frac{\rho}{2} Q_2,
		\end{aligned}
		\right.
	\end{equation}
	%\begin{equation}
	%\begin{aligned}
	where $\sigma_{tr,g}=\frac{1}{\Delta E_g}\int_{E_{g-\frac{1}{2}}}^{E_{g+\frac{1}{2}}}  \sigma_{tr} \text{d}E$ and the momentum transfer cross section $\sigma_{tr}$ is defined in the appendix.

	%%%%%%%%%%%%%%%%%%%%%%%%%%%%%%%%%%%%%%%%%%%%%%%%%%%%%%%%%%
	\subsection{The depth splitting method}
	All protons enter at $x=0$ and the computational domain in the $x$ variable is $[0,L]$, the advection coefficient in front of $\partial_x$ in \eqref{eqn:0041-d} is a constant $1$, thus the $x$ variable is similar to the time variable in classical time evolutionary problems and the injected proton distribution can be considered as the initial condition.  We rewrite  \eqref{eqn:0041-d} into the following operator form:
	\begin{equation}\label{eqn:0041-do}
		\left\{
		\begin{aligned}
			&\frac{\partial\psi_g^1}{\partial x}= \mathcal{L}_1\psi_g^1+\sigma_{tr,g}\rho\mathcal{L}_2\psi_g^1+\mathcal{L}_{3,g}(\psi_g^1,\psi_g^2,\psi_{g\pm 1}^1,\psi_{g\pm 1}^2), \\
			&\frac{\partial\psi_g^2}{\partial x}= \mathcal{L}_1\psi_g^2+\sigma_{tr,g}\rho\mathcal{L}_2\psi_g^2+\mathcal{L}_{4,g}(\psi_g^1,\psi_g^2,\psi_{g\pm 1}^1,\psi_{g\pm 1}^2),
		\end{aligned}
		\right.
	\end{equation}
	with
	$$
	\mathcal{L}_1=-u\partial_y-v\partial_z,\quad \mathcal{L}_2=\frac{1}{2}(\partial_u^2+\partial_v^2)
	$$and
	$$\begin{aligned}
		\mathcal{L}_{3,g}(\psi_g^1,\psi_g^2,\psi_{g\pm 1}^1,\psi_{g\pm 1}^2)=&\rho\frac{S_{g+\frac{1}{2}}}{\Delta E_g}\left(\psi_{ g+1}^1-\psi_{ g+1}^2\right)-\rho\frac{S_{g-\frac{1}{2}}}{\Delta E_g}\left(\psi_{ g}^1-\psi_{ g}^2\right)-\sigma_{c,t,g}\psi^1_g\\
		&+\frac{\rho}{2}\frac{T_{g+1}-T_g}{(\Delta E_g)^2}\left(\psi_{ g+1}^1-\psi_{ g+1}^2\right)-\frac{\rho}{2}\frac{T_{g}-T_{g-1}}{(\Delta E_g)^2}\left(\psi_{ g}^1-\psi_{ g}^2\right)+\frac{\rho}{2} Q_1,\\
		\mathcal{L}_{4,g}(\psi_g^1,\psi_g^2,\psi_{g\pm 1}^1,\psi_{g\pm 1}^2)=&3 \rho\frac{S_{g+\frac{1}{2}}}{\Delta E_g}\left(\psi^1_{g+1}-\psi^2_{g+1}\right)
		-3\rho \frac{S_{g-\frac{1}{2}}}{\Delta E_g}\left(\psi^1_{g}-\psi^2_{g}\right)-3 \rho\frac{S_{g+\frac{1}{2}}+S_{g-\frac{1}{2}}}{\Delta E_g} \psi^1_{g}\\
		&-\rho\frac{S_{g+\frac{1}{2}}-S_{g-\frac{1}{2}}}{\Delta E_g} \psi^2_{g}-\sigma_{c,t,g}\psi^2_g
		+\frac{3\rho}{2} \frac{T_{g+1}-T_{g}}{(\Delta E_g)^2}\left(\psi^1_{g+1}-\psi^2_{g+1}\right)\\
		&-\frac{3\rho}{2} \frac{T_{g}-T_{g-1}}{(\Delta E_g)^2}\left(\psi^1_{g}-\psi^2_{g}\right)-\frac{3\rho}{2} \frac{T_{g+1}-T_{g-1}}{(\Delta E_g)^2}\psi^1_{g}+\frac{\rho}{2} Q_2.
	\end{aligned}
	$$
	One can then discretize the system \eqref{eqn:0041-do} by the classical  splitting method. More precisely, let $[0,L]$ be divided into $N_x$ intervals and the size of each interval is $\Delta x=L/N_x$. Let $\psi_g^{l,s}(y,z,u,v)\approx\psi_g^l(s\Delta x,y,z,u,v)$, ($l=1,2$). In order to get $\psi_g^{l,s+1}$ starting from $\psi_g^{l,s}$,
	the first order splitting method writes:
	
	Let $\psi_g^{l,s(1)}(s\Delta x)=\psi_g^{l,s}$, solve
	\begin{equation}\label{split-1}
		\left\{
		\begin{aligned}
			\frac{\partial\psi_g^{1,s(1)}}{\partial x}&=\mathcal{L}_3(\psi_g^{1,s(1)},\psi_g^{2,s(1)},\psi_{g\pm 1}^{1,s(1)},\psi_{g\pm 1}^{2,s(1)}),\\
			\frac{\partial\psi_g^{2,s(1)}}{\partial x}&=\mathcal{L}_4(\psi_g^{1,s(1)},\psi_g^{2,s(1)},\psi_{g\pm 1}^{1,s(1)},\psi_{g\pm 1}^{2,s(1)}),
		\end{aligned}
		\right.
	\end{equation}
	for one step from $s\Delta x$ to $(s+1)\Delta x$.
	After getting $\psi_g^{l,s+1(1)}\approx\psi_g^l\big((s+1)\Delta x,y,z,u,v\big)$, one updates the second system for one step from $s\Delta x$ to $(s+1)\Delta x$, starting from $\psi_g^{l,s+1(1)}$:
	\begin{equation}\label{split-2}
		\left\{
		\begin{aligned}
			&\frac{\partial\psi_g^{1,s(2)}}{\partial x}=\mathcal{L}_1\psi_g^{1,s(2)},\\
			&\frac{\partial \psi_g^{1,s(2)}}{\partial x}=\mathcal{L}_1\psi_g^{2,s(2)}.
		\end{aligned}
		\right.
	\end{equation}
	Finally, using $\psi_g^{l,s+1(2)}\approx\psi_g^{l,s(2)}\big((s+1)\Delta x,y,z,u,v\big)$ as the initial condition, one updates the third system for one step, from $s\Delta x$ to $(s+1)\Delta x$:
	\begin{equation}\label{split-3}
		\left\{
		\begin{aligned}
			&\frac{\partial\psi_g^{1,s(3)}}{\partial x}=\sigma_{tr,g}\rho\mathcal{L}_2\psi_g^{1,s(3)},\\
			&\frac{\partial \psi_g^{2,s(3)}}{\partial x}=\sigma_{tr,g}\rho\mathcal{L}_2\psi_g^{2,s(3)}.
		\end{aligned}
		\right.
	\end{equation}
	One can get $\psi_g^{l,s+1}=\psi_g^{l,s+1(3)}$.

	%
	%
	%\begin{subequations}\label{eqn:004}
	%	\begin{numcases}{}
		%	\frac{\psi^{s+1}-\psi^{s}}{\Delta x}=\rho\frac{\partial}{\partial E}\left(S(\textbf{r},E)\psi\right)+\rho\frac{1}{2}\frac{\partial^2}{\partial E^2}\left(T(\textbf{r},E)\psi\right), \label{eqn:case11004}\\
		%	\frac{\psi^{s+1}-\psi^{s}}{\Delta x}+u\frac{\partial\psi}{\partial y}+v\frac{\partial\psi}{\partial z}=0,\label{eqn:case12004}\\
		%	\frac{\psi^{s+1}-\psi^{s}}{\Delta x}=\rho\frac{\sigma_{tr}}{2}\left(\frac{\partial^2}{\partial u^2}+\frac{\partial^2}{\partial v^2}\right)\psi,\label{eqn:case13004}
		%	\end{numcases}
	%\end{subequations}
In order to achieve second order convergence	in the $x$ variable, one can use second order Strang splitting method as in \cite{spiteri2023beyond} and Crank-Nicolson (C-N) method for each subsystem. More precisely, starting from $\psi_g^{l,s}$ ($l=1,2$), one solve \eqref{split-1} for half step $\Delta x/2$ by CN method:
	\begin{equation*}
		\left\{
		\begin{aligned}
			&\frac{\bar\psi_g^{1,s(1)}-\psi_g^{1,s}}{\Delta x/2}=\frac{1}{2}\left(\mathcal{L}_3(\bar\psi_g^{1,s(1)},\bar\psi_g^{2,s(1)},\bar \psi_{g\pm 1}^{1,s(1)},\bar \psi_{g\pm 1}^{2,s(1)})+\mathcal{L}_3(\psi_g^{1,s},\psi_g^{2,s},\psi_{g\pm 1}^{1,s},\psi_{g\pm 1}^{2,s})\right),\\
			&\frac{\bar\psi_g^{2,s(1)}-\psi_g^{2,s}}{\Delta x/2}=\frac{1}{2}\left(\mathcal{L}_4(\bar\psi_g^{1,s(1)},\bar\psi_g^{2,s(1)},\bar \psi_{g\pm 1}^{1,s(1)},\bar \psi_{g\pm 1}^{2,s(1)})+\mathcal{L}_4(\psi_g^{1,s},\psi_g^{2,s},\psi_{g\pm 1}^{1,s},\psi_{g\pm 1}^{2,s})\right).
		\end{aligned}
		\right.
	\end{equation*}
	Then, after getting $(\bar\psi_g^{1,s(1)},\bar\psi_g^{2,s(1)})$, one updates the second system for half step:
	\begin{equation*}
		\left\{
		\begin{aligned}
			&\frac{\bar\psi_g^{1,s(2)}-\bar\psi_g^{1,s(1)}}{\Delta x/2}=\frac{1}{2}\left(\mathcal{L}_1\bar\psi_g^{1,s(2)}+\mathcal{L}_1\bar\psi_g^{1,s(1)}\right),\\
			&\frac{\bar\psi_g^{2,s(2)}-\bar\psi_g^{2,s(1)}}{\Delta x/2}=\frac{1}{2}\left(\mathcal{L}_1\bar\psi_g^{2,s(2)}+\mathcal{L}_1\bar\psi_g^{2,s(1)}\right).
		\end{aligned}
		\right.
	\end{equation*}
	Using $(\bar\psi_g^{1,s(2)},\bar\psi_g^{2,s(2)})$, one updates the third system for one step:
	\begin{equation*}
		\left\{
		\begin{aligned}
			&\frac{\bar\psi_g^{1,s(3)}-\bar\psi_g^{1,s(2)}}{\Delta x}=\frac{\sigma_{tr,g}\rho}{2}\left(\mathcal{L}_2\bar\psi_g^{1,s(3)}+\mathcal{L}_2\bar\psi_g^{1,s(2)}\right),\\
			&\frac{\bar\psi_g^{2,s(3)}-\bar\psi_g^{2,s(2)}}{\Delta x}=\frac{\sigma_{tr,g}\rho}{2}\left(\mathcal{L}_2\bar\psi_g^{2,s(3)}+\mathcal{L}_2\bar\psi_g^{2,s(2)}\right).
		\end{aligned}
		\right.
	\end{equation*}
	One updates the forth system for half step:
	\begin{equation*}
		\left\{
		\begin{aligned}
			&\frac{\bar\psi_g^{1,s(4)}-\bar\psi_g^{1,s(3)}}{\Delta x/2}=\frac{1}{2}\left(\mathcal{L}_1\bar\psi_g^{1,s(4)}+\mathcal{L}_1\bar\psi_g^{1,s(3)}\right),\\
			&\frac{\bar\psi_g^{2,s(4)}-\bar\psi_g^{2,s(3)}}{\Delta x/2}=\frac{1}{2}\left(\mathcal{L}_1\bar\psi_g^{2,s(4)}+\mathcal{L}_1\bar\psi_g^{2,s(3)}\right).
		\end{aligned}
		\right.
	\end{equation*}
	At last, one updates the fifth system for half step:
	\begin{equation*}
		\left\{
		\begin{aligned}
			&\frac{\bar\psi_g^{1,s(5)}-\bar\psi_g^{1,s(4)}}{\Delta x/2}=\frac{1}{2}\left(\mathcal{L}_3(\bar\psi_g^{1,s(5)},\bar\psi_g^{2,s(5)},\bar\psi_{g\pm 1}^{1,s(5)},\bar\psi_{g\pm 1}^{2,s(5)})+\mathcal{L}_3(\bar\psi_g^{1,s(4)},\bar\psi_g^{2,s(4)},\bar\psi_{g\pm 1}^{1,s(4)},\bar\psi_{g\pm 1}^{2,s(4)})\right),\\
			&\frac{\bar\psi_g^{2,s(5)}-\bar\psi_g^{2,s(4)}}{\Delta x/2}=\frac{1}{2}\left(\mathcal{L}_4(\bar\psi_g^{1,s(5)},\bar\psi_g^{2,s(5)},\bar\psi_{g\pm 1}^{1,s(5)},\bar\psi_{g\pm 1}^{2,s(5)})+\mathcal{L}_4(\bar\psi_g^{1,s(4)},\bar\psi_g^{2,s(4)},\bar\psi_{g\pm 1}^{1,s(4)},\bar\psi_{g\pm 1}^{2,s(4)})\right).
		\end{aligned}
		\right.
	\end{equation*}
	Finally, one has $\psi_g^{l,s+1}=\bar\psi_g^{l,s(5)}$.
	%
	%
	%why you can consider $x=0$ as the initial condition
	%\textcolor{red}{The injection of proton is from $x=0$, the proton can have a distribution in $E$, $u$,$v$, what does it mean for pronton to have a distribution, why you can consider $x=0$ as the initial condition? approximation is important, inflow is usually concentrate at some point, spot of the machine.}
	%and by using the operator splitting method to split the system. 
	%\textcolor{red}{the discretization in $x$, learn from time splitting method, strang time splitting method.}
	%Then one can get the following system
\begin{remark}	
In equation \eqref{eqn:pri-se2}, there is a source term about the primary proton 
$$
\int_{4\pi}\int_E^\infty\sigma_{c,s}\psi_p\text{d}u'\text{d}v'\text{d}E',
$$
and it can be added to \eqref{split-1}. We will introduce the energy discretization of the source term. Firstly, testing the source term by $\frac{\phi_g^1}{\Delta E_g}$, one can get
	$$\begin{aligned}
			&\frac{1}{\Delta E_g} \int_{E_{g-\frac{1}{2}}}^{E_{g+\frac{1}{2}}} \phi_g^1  \int_{4\pi}\int_E^\infty\sigma_{c,s}\psi\text{d}u'\text{d}v'\text{d}E'\text{d}E\\
   =&\frac{1}{\Delta E_g} \int_{E_{g-\frac{1}{2}}}^{E_{g+\frac{1}{2}}} \phi_g^1  \Delta u\Delta v\sum_{u,v}\sum_{g'=G}^{g}\sigma_{c,s,g'\to g}\psi^1_{g'}\Delta E_{g'}\text{d}E\\  
   =&\Delta u\Delta v\sum_{u,v}\sum_{g'=G}^{g}\sigma_{c,s,g'\to g}\psi^1_{g'}\frac{\Delta E_{g'}}{\Delta E_{g}},
  % \approx &\frac{2}{\Delta E_g}\frac{T_{g+1} +T_{g} }{2} \frac{ \psi^1_{g+1}- \psi^1_g}{\Delta E_g+\Delta E_{g+1}}-\frac{2}{\Delta E_g}\frac{T_{g-1} +T_{g} }{2} \frac{\psi^1_g-\psi^1_{g-1}}{\Delta E_g+\Delta E_{g-1}},
		\end{aligned}$$
and test the source term by $\frac{3\phi_g^2}{\Delta E_g}$ which yields
$$\begin{aligned}
			&\frac{1}{\Delta E_g} \int_{E_{g-\frac{1}{2}}}^{E_{g+\frac{1}{2}}} \phi_g^2  \int_{4\pi}\int_E^\infty\sigma_{c,s}\psi\text{d}u'\text{d}v'\text{d}E'\text{d}E
   \approx 0,
		\end{aligned}$$
where we assume that  $$\int_{E_{g-\frac{1}{2}}}^{E_{g+\frac{1}{2}}} \phi_g^2  \sigma_{c,s}\text{d}E=0.$$
\end{remark}
\begin{remark}	
In the origin system, the variables: energy, space and angular are coupled together, while the depth splitting system is easy to parallel. For the subsystem \eqref{split-1}, suppose the angular and YZ plane has $N_u\times N_v\times N_y\times N_z$ different grids,  \eqref{split-1} can be solved in parallel for each grid point. Similarly, the subsystem \eqref{split-2} can be solved in parallel for the $N_u\times N_v\times G$ different grids on  the energy and angular plane.  Subsystem \eqref{split-3} can be solved in parallel for the $N_y\times N_z\times G$ different grids on the energy and $YZ$ plane.
  \end{remark}
% \begin{remark}
%	From Fig.~\ref{fig:s_c}, we can see that the straggling coefficient $T$ is small. Moreover, the straggling operator acts as a diffusion operator  which increases the computation time significantly. Thus, we treat the straggling operator explicitly in the implicit part of the CN method.   
%\end{remark}
 %   \textcolor{red}{For most types of radiation treatment, the Fermi pencil-beam model  offers a   fast and reliable alternative. Despite its success in most clinical uses, it fails in complicated settings like inhomogeneities. The Fermi–Eyges approximation allows only for the case $\sigma_{tr}=\sigma_{tr}(x)$, i.e., for layered heterogeneities. Its application  to nonlayered media is difficult. However, we use numerical scheme to treat $\sigma_{tr}$, which is not involved in Fermi–Eyges  theory. Thus our scheme can handle the case $\sigma_{tr}=\sigma_{tr}(x,y,z)$.}
  
%\textcolor{red}{make this paragraph a remark, find out why this commend does not work}

\subsection{Angular and $YZ$ discretization}
Since zero boundary condition in the angular variable $(u,v)$ is used, we solve the subsystem \eqref{split-3} by using the Fourier based methods.

Firstly, we use the finite difference method to discretize the equation. Let the computational domain in the $(u,v)$ variables be $[a,b]\times[c,d]$. $N_u$ and $N_v$ grid points are used for $u,v$  respectively and the mesh sizes are $\Delta u=\frac{b-a}{N_u}$ and $\Delta v=\frac{d-c}{N_v}$. 
%\begin{subequations}\label{eqn:0041-ds}
%\begin{numcases}{}
%\frac{\partial\psi_{l,j,g}^1}{\partial x}=\rho\frac{S_{l,j,g+\frac{1}{2}}}{\Delta E_g}\left(\psi_{ l,j,g+1}^1-\psi_{ l,j,g+1}^2\right)-\rho\frac{S_{l,j,g-\frac{1}{2}}}{\Delta E_g}\left(\psi_{ l,j,g}^1-\psi_{ l,j,g}^2\right)\nonumber\\
%\hspace{1.2cm}+\frac{\rho}{\Delta E_g} \frac{\left(T \psi^1\right)_{l,j,g+1}-\left(T \psi^1\right)_{l,j,g}}{\Delta E_g+\Delta E_{g+1}}-\frac{\rho}{\Delta E_g} \frac{\left(T \psi^1\right)_{l,j,g}-\left(T \psi^1\right)_{l,j,g-1}}{\Delta E_g+\Delta E_{g-1}}, \label{eqn:case110041-ds}\\
%\frac{\partial\psi_{l,j,g}^2}{\partial x}=3 \frac{S_{l,j,g+\frac{1}{2}}}{\Delta E_g}\left(\psi^1_{l,j,g+1}-\psi^2_{l,j,g+1}\right)
%-3 \frac{S_{l,j,g-\frac{1}{2}}}{\Delta E_g}\left(\psi^1_{l,j,g}-\psi^2_{l,j,g}\right)\nonumber\\
%\hspace{1.2cm}-3 \frac{S_{l,j,g+\frac{1}{2}}+S_{l,j,g-\frac{1}{2}}}{\Delta E_g} \psi^1_{l,j,g}-\frac{S_{l,j,g+\frac{1}{2}}-S_{l,j,g-\frac{1}{2}}}{\Delta E_g} \psi^2_{l,j,g}.\label{eqn:case120041-ds}\\
%\frac{\partial\psi_{l,j,g}^1}{\partial x}+u_l\frac{\partial\psi_{l,j,g}^1}{\partial y}+v_j\frac{\partial\psi_{l,j,g}^1}{\partial z}=0,\label{eqn:case130041-ds}\\
%\frac{\partial\psi_{l,j,g}^2}{\partial x}+u_l\frac{\partial\psi_{l,j,g}^2}{\partial y}+v_j\frac{\partial\psi_{l,j,g}^2}{\partial z}=0.\label{eqn:case140041-ds}
% \end{numcases}
%\end{subequations}
By using central difference, the discretization of \eqref{split-3} is as follow:
\begin{equation}\label{eqn:sd}
\left\{
    \begin{aligned}
   &  \frac{\partial\psi_{i,j,g}^{1,s(3)}}{\partial x}=\frac{\sigma_{tr,g}\rho}{2}\left(\frac{\psi_{i+1,j,g}^{1,s(3)}-2\psi_{i,j,g}^{1,s(3)}+\psi_{i-1,j,g}^{1,s(3)}}{(\Delta u)^2} +\frac{\psi_{i,j+1,g}^{1,s(3)}-2\psi_{i,j,g}^{1,s(3)}+\psi_{i,j-1,g}^{1,s(3)}}{(\Delta v)^2}\right),\\
    & \frac{\partial\psi_{i,j,g}^{2,s(3)}}{\partial x}=\frac{\sigma_{tr,g}\rho}{2}\left(\frac{\psi_{i+1,j,g}^{2,s(3)}-2\psi_{i,j,g}^{2,s(3)}+\psi_{i-1,j,g}^{2,s(3)}}{(\Delta u)^2} +\frac{\psi_{i,j+1,g}^{2,s(3)}-2\psi_{i,j,g}^{2,s(3)}+\psi_{i,j-1,g}^{2,s(3)}}{(\Delta v)^2}\right),
     \end{aligned}
    \right.
\end{equation}
where $\psi_{i,j,g}^{l,s(3)}$ ($l=1,2$) is the approximation  of $\psi_{g}^{l,s(3)}(u_i,v_j)$.
We will illustrate how to solve the above   equations by using the  Fourier based method. Due to the zero boundary condition, we utilize the discrete Sine series of $\{\psi_{i,j,g}^{l,s(3)}\}$ ($l=1,2$) as follows:
\begin{equation}\label{eqn_fly}
    \psi_{i,j,g}^{l,s(3)}=\frac{1}{IJ}\sum_{m=0}^{I-1}\sum_{n=0}^{J-1}\hat{\psi}_{m,n,g}^{l,s(3)}\text{sin}\frac{\pi im}{I}\text{sin}\frac{\pi jn}{J}.
\end{equation}
Substituting \eqref{eqn_fly} into \eqref{eqn:sd}, one can obtain
\begin{equation*}
\left\{
    \begin{aligned}
   &  \frac{\partial\hat{\psi}_{m,n,g}^{1,s(3)}}{\partial x}=\frac{\sigma_{tr,g}\rho}{2}\left[\frac{\hat{\psi}_{m,n,g}^{1,s(3)}}{(\Delta u)^2}\left(2\text{cos}\frac{\pi m}{I}-2\right)
  +\frac{\hat{\psi}_{m,n,g}^{1,s(3)}}{(\Delta v)^2}\left(2\text{cos}\frac{\pi n}{J}-2\right)\right],\\
  & \frac{\partial\hat{\psi}_{m,n,g}^{2,s(3)}}{\partial x}=\frac{\sigma_{tr,g}\rho}{2}\left[\frac{\hat{\psi}_{m,n,g}^{2,s(3)}}{(\Delta u)^2}\left(2\text{cos}\frac{\pi m}{I}-2\right)
  +\frac{\hat{\psi}_{m,n,g}^{2,s(3)}}{(\Delta v)^2}\left(2\text{cos}\frac{\pi n}{J}-2\right)\right].
    \end{aligned}
    \right.
    \end{equation*}
Then one can get $\hat{\psi}_{m,n,g}^{l,s(3)}$  for the next depth by an ODE solver, and $\psi_{m,n,g}^{l,s(3)}$  can be obtained by inverse transform. Moreover, the transform and inverse transform can use the fast Fourier method.
    
%\subsection{Y-Z discretization}
The discretization  in $YZ$ is achieved by using a second order finite volume scheme. We only focus on the discretization of \eqref{split-2}. Let $(y_{p+1/2}, z_{q+1/2})$ be the grid points of the $YZ$  domain. For uniform mesh $\Delta y=y_{p+1/2}-y_{p-1/2}$ and $\Delta z=z_{q+1/2}-z_{q-1/2}$, let $\psi_{i,j,p,q,g}^{l,s(2)}=\frac{1}{\Delta y\Delta z}\int_{y_{p-1/2}}^{y_{p+1/2}}\int_{z_{q-1/2}}^{z_{q+1/2}}\psi_{i,j,g}^{l,s(2)}\text{d}z\text{d}y$, the second order MUSCL scheme writes \cite{kurganov2000new}:
%\begin{equation*}
%\left\{
%    \begin{aligned}
%    \int_{y_{p-1/2}}^{y_{p+1/2}}\int_{z_{q-1/2}}^{z_{q+1/2}}\frac{\partial\psi_{i,j,g}^{1,s(2)}}{\partial x}\text{d}z\text{d}y+\int_{y_{p-1/2}}^{y_{p+1/2}}\int_{z_{q-1/2}}^{z_{q+1/2}}u_i\frac{\partial\psi_{i,j,g}^{1,s(2)}}{\partial y}\text{d}z\text{d}y+\int_{y_{p-1/2}}^{y_{p+1/2}}\int_{z_{q-1/2}}^{z_{q+1/2}}v_j\frac{\partial\psi_{i,j,g}^{1,s(2)}}{\partial z}\text{d}z\text{d}y=0,\\
 %   \int_{y_{p-1/2}}^{y_{p+1/2}}\int_{z_{q-1/2}}^{z_{q+1/2}}\frac{\partial\psi_{i,j,g}^{2,s(2)}}{\partial x}\text{d}z\text{d}y+\int_{y_{p-1/2}}^{y_{p+1/2}}\int_{z_{q-1/2}}^{z_{q+1/2}}u_i\frac{\partial\psi_{i,j,g}^{2,s(2)}}{\partial y}\text{d}z\text{d}y+\int_{y_{p-1/2}}^{y_{p+1/2}}\int_{z_{q-1/2}}^{z_{q+1/2}}v_j\frac{\partial\psi_{i,j,g}^{2,s(2)}}{\partial z}\text{d}z\text{d}y=0,
%    \end{aligned}
%    \right.
%\end{equation*}

\begin{equation}\label{distansport}
	\left\{
	\begin{aligned}
	\frac{\partial\psi_{i,j,p,q,g}^{1,s(2)}}{\partial x}+\frac{\Psi^{1,s(2)}_{i,j,p+1/2,q,g}-\Psi^{1,s(2)}_{i,j,p-1/2,q,g}}{\Delta y}+\frac{\Psi^{1,s(2)}_{i,j,p,q+1/2,g}-\Psi^{1,s(2)}_{i,j,p,q-1/2,g}}{\Delta z}=0,\\
		\frac{\partial\psi_{i,j,p,q,g}^{2,s(2)}}{\partial x}+\frac{\Psi^{2,s(2)}_{i,j,p+1/2,q,g}-\Psi^{2,s(2)}_{i,j,p-1/2,q,g}}{\Delta y}+\frac{\Psi^{2,s(2)}_{i,j,p,q+1/2,g}-\Psi^{2,s(2)}_{i,j,p,q-1/2,g}}{\Delta z}=0,
	\end{aligned}
	\right.
\end{equation}
where the fluxes are defined as
$$
\Psi^{l,s(2)}_{i,j,p+1/2,q,g}=u^{+}_i\psi^{l,s(2),L}_{i,j,p+1/2,q,g}-u^{-}_i\psi^{l,s(2),R}_{i,j,p+1/2,q,g},\quad \Psi^{l,s(2)}_{i,j,p,q+1/2,g}=v^{+}_j\psi^{l,s(2),L}_{i,j,p,q+1/2,g}-v^{-}_j\psi^{l,s(2),R}_{i,j,p,q+1/2,g},
$$
with $u^{+}_i=\text{max}(u_i,0)$, $u^{-}_i=\text{max}(-u_i,0)$, $v^{+}_j=\text{max}(v_j,0)$, $v^{-}_j=\text{max}(-v_j,0)$, and
$$
\begin{aligned}
\psi^{l,s(2),L}_{i,j,p+1/2,q,g}=\psi^{l,s(2)}_{i,j,p,q,g}+\frac{\delta^{l,s(2)}_{i,j,p+1/2,q,g}}{2},\quad \psi^{l,s(2),R}_{i,j,p+1/2,q,g}=\psi^{l,s(2)}_{i,j,p+1,q,g}-\frac{\delta^{l,s(2)}_{i,j,p+1/2,q,g}}{2},\\
\psi^{l,s(2),L}_{i,j,p,q+1/2,g}=\psi^{l,s(2)}_{i,j,p,q,g}+\frac{\delta^{l,s(2)}_{i,j,p,q+1/2,g}}{2},\quad \psi^{l,s(2),R}_{i,j,p,q+1/2,g}=\psi^{l,s(2)}_{i,j,p+1,q,g}-\frac{\delta^{l,s(2)}_{i,j,p,q+1/2,g}}{2}.
\end{aligned}
$$
Here $\delta^{l,s(2)}$ is a slope with a limiter, as in \cite{leveque2002finite}, which is defined by
$$
\delta^{l,s(2)}_{i,j,p+1/2,q,g}=(\psi^{l,s(2)}_{i,j,p+1,q,g}-\psi^{l,s(2)}_{i,j,p,q,g})\eta(\theta^{l,s(2)}_{i,j,p,q,g}),
$$
with 
$$
\eta(\theta^{l,s(2)}_{i,j,p,q,g})=\text{max}[0,\text{min}(2\theta^{l,s(2)}_{i,j,p,q,g},1),\text{min}(\theta^{l,s(2)}_{i,j,p,q,g},2)],\quad \theta^{l,s(2)}_{i,j,p,q,g}=\frac{\psi^{l,s(2)}_{i,j,p,q,g}-\psi^{l,s(2)}_{i,j,p-1,q,g}}{\psi^{l,s(2)}_{i,j,p+1,q,g}-\psi^{l,s(2)}_{i,j,p,q,g}}.
$$
\begin{remark}
	  The CN depth discretization can be employed in both \eqref{eqn:sd} and \eqref{distansport}. Based on the Fourier based method, the implicit part of the CN depth discretization can be solved explicitly. Moreover, one can treat the slope explicitly in the implicit part of the CN method for \eqref{distansport}, which allows for an explicit treatment of the subsystem for spatial transportation.
	\end{remark}

\section{Numerical Tests}
%\textcolor{red}{include information on the number of spatial "steps" required in their calculations.}
%\textcolor{red}{Are the beams the authors are using of a diameter appropriate for proton therapy}
We conduct numerical tests for the performance of our proposed method. The reference solution is calculated by FLUKA. FLUKA is a general Monte Carlo code designed for particle transportation and interactions with background material. Initially developed for accelerator shielding, FLUKA has been expanded to a wide range of applications, including cosmic rays and radiation therapy. FLUKA can provide directly the dose distribution of the  protons \cite{battistoni2007fluka}. %\textcolor{blue}{Moreover, we compare the dose obtained by the model in \eqref{f-e} and the semi-analytical approximation derived in \cite{burlacu2023deterministic}.}

As in \cite{burlacu2023deterministic,embriaco2018monet}, we will consider the energy deposition of a single Gaussian beam with various energy distributions and being injected into different materials.  In all numerical examples, the beam enters the target at $(x,y,z)=(0,0,0)$; the computational domain is $(x,y,z)\in  (0\ \text{cm},40\ \text{cm})\times (-4\ \text{cm},4\ \text{cm})\times(-4\ \text{cm},4\ \text{cm})$. The domain can be composed of different materials. Moreover, since the protons have highly forward peaked scattering, for the angular variables, we use $(u,v)\in (-1,1)\times(-1,1)$. The energy domain is set to be [1 MeV, 260 MeV]. On the left boundary, the boundary condition is modeled by a Gaussian distribution:
\begin{equation}\label{bdy}
    \psi(0,y,z,u,v,E)=Ce^{-\frac{y^2}{2\sigma_y^2}-\frac{u^2}{2\sigma_u^2}}e^{-\frac{z^2}{2\sigma_z^2}-\frac{v^2}{2\sigma_v^2}}e^{-\frac{(E-E_0)^2}{2\sigma_E^2}},
\end{equation}
where $C$ is a constant, $\sigma_y=\sigma_z=0.3\ \text{cm}$, $\sigma_u=\sigma_v=10^{-6}$, $\sigma_E=1\ \text{MeV}$, and $E_0$ is determined by the energy of the incoming beam. Three different beams $E_0=50 \text{MeV},100 \text{MeV},230 \text{MeV}$ and three different materials are tested in the subsequent part. The characteristics of water, bone and air are listed in Table.~\ref{tab:stab2}.
\begin{table}[!ht]
	\centering
	\caption{Characteristics of material target. }\label{tab:stab2}	
	\begin{tabular}{lll}
				\hline
				Target &density (\text{g cm$^{-3}$})&composition (weight fraction)\\
				\hline
				water& 1.000&H 0.1111, O 0.8889\\
\hline
			bone& 1.757& \tabincell{cc}{H 0.04200, C 0.1940, N 0.04000,\\
				O 0.4250, Na 0.001000, Mg 0.002000,\\
				P 0.09200, S 0.003000, Ca 0.2010 }\\
				\hline
			air& 0.001205& C 0.0001248, N 0.7553, O 0.2318, Ar 0.01283\\
			\hline
			\end{tabular}\vspace{0cm}
			\end{table}

For each energy group, we have that
\begin{equation*}
    \sigma_{c,s,g}=\sigma_{c,t,g}\mathcal{P}\left(g'\to g,(u',v')\to(u,v)\right),
\end{equation*}
where $\sigma_{c,s,g}$ and $\sigma_{c,t,g}$ are respectively the catastrophic scattering and total cross sections; $\mathcal{P}$ is the probability transition matrix. One can assume that
\begin{equation*}
    \mathcal{P}\left(g'\to g,(u',v')\to(u,v)\right)=\mathcal{P}^1\left(g'\to g\right)\mathcal{P}^2_g\left((u',v')\to(u,v)\right),
\end{equation*}
where $\mathcal{P}^1\left(g'\to g\right)$ is the probability that the protons in energy group $g'$ are scattered to the energy group $g$, and $\mathcal{P}^2_g\left((u',v')\to(u,v)\right)$ is the probability that the protons in direction $(u',v')$ are scattered to the direction $(u,v)$ in the energy group $g$. We obtain the cross sections by collecting the statistical data in FLUKA. The detail can be found in Appendix.~\ref{sec-d}.

	%\subsection{Dose}   

			The role of dose in radiotherapy is crucial, as it directly impacts the effectiveness of the treatment and the side effects. The dose of proton therapy refers to the amount of radiation delivered to or energy deposited in a specific area. The dose is defined by
			$$
	    D(\textbf{r})=\frac{E_{\text{dep}}(\textbf{r})}{\rho(\textbf{r})},
			$$
			where $E_{\text{dep}}$ is the energy  deposition  in the specific region, which includes primary proton energy  deposition and secondary proton energy  deposition. Here we present the formula for the energy deposition of primary protons. A similar formula for secondary protons can be obtained.
   %which is related to the energy flux 
%			$$
%			\int_{-\infty}^{+\infty}\int_{-\infty}^{+\infty}\int_{E_{min}}^{E_{max}}E\psi\text{d}E\text{d}u\text{d}v.
%			$$ 
%			The energy flux there  can be gotten 
Multiplying  equation \eqref{f-e} by $E$, and integrating over the energy and angular variables, one can get 
			\begin{equation*}
				\begin{aligned}
					&\int_{-\infty}^{+\infty}\int_{-\infty}^{+\infty}\int_{E_{min}}^{E_{max}}E\frac{\partial\psi_p}{\partial x}\text{d}E\text{d}u\text{d}v\\
					=&\int_{-\infty}^{+\infty}\int_{-\infty}^{+\infty}\int_{E_{min}}^{E_{max}}\rho E\frac{\partial}{\partial E}\left(S(\textbf{r},E)\psi_p\right)+\rho E\frac{1}{2}\frac{\partial^2}{\partial E^2}\left(T(\textbf{r},E)\psi_p\right)-E\sigma_{c,t}\psi_p \text{d}E\text{d}u\text{d}v.
				\end{aligned}
			\end{equation*}
		The energy deposition of primary protons $E_{\text{dep},p}$ is the difference of the energy flux between two successive depth, which is given by $-\int_{-\infty}^{+\infty}\int_{-\infty}^{+\infty}\int_{E_{min}}^{E_{max}}E\frac{\partial\psi_p}{\partial x}\text{d}E\text{d}u\text{d}v$, thus 
			\begin{equation}\label{en_d}
			E_{\text{dep},p}(\textbf{r})=-\int_{-\infty}^{+\infty}\int_{-\infty}^{+\infty}\int_{E_{min}}^{E_{max}}\rho E\frac{\partial}{\partial E}\left(S(\textbf{r},E)\psi\right)+\rho E\frac{1}{2}\frac{\partial^2}{\partial E^2}\left(T(\textbf{r},E)\psi\right)-E\sigma_{c,t}\psi_p \text{d}E\text{d}u\text{d}v.
		\end{equation}

 In the subsequent part, we will compare the integrated depth dose (IDD), spot distribution and the dose distribution in the $XY$ plane. 
For the spatial variables $(x,y,z)$, we use a Cartesian $4000 \times 80 \times 80$ mesh. The energy domain $[1 \text{MeV}, 260 \text{MeV}]$ is divided into $500$ groups. For the angular variable $(u,v)$, $20\times 20$ uniform Cartesian mesh is used. These choices are to show the agreement of with FLUKA, the meshes can be coarser if the efficiency is more important.

\subsection{Test of Convergence}
\begin{figure}[!h]
	\centering
	\includegraphics[width=0.6\textwidth]{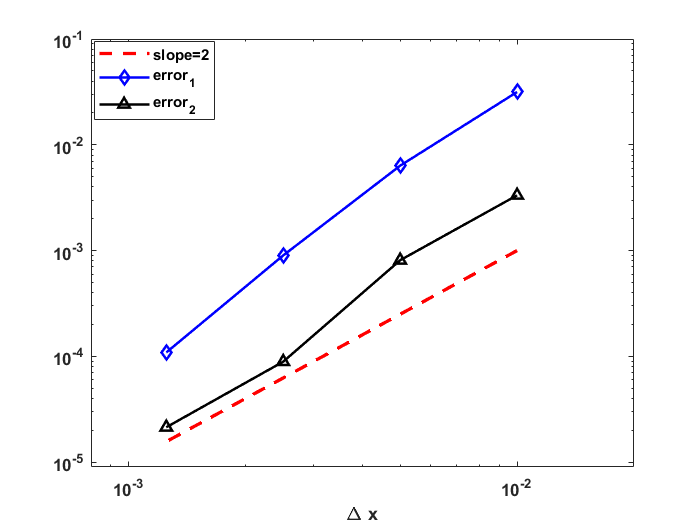}
	\caption{Plot of error \eqref{error} with decreasing $\Delta x$ and $\Delta E_g$.}
	\label{fig:conv}
\end{figure}

\begin{figure}[!h]
\centering
 \subfloat   % 
  {
      \includegraphics[width=2.15in, height=2.5in]{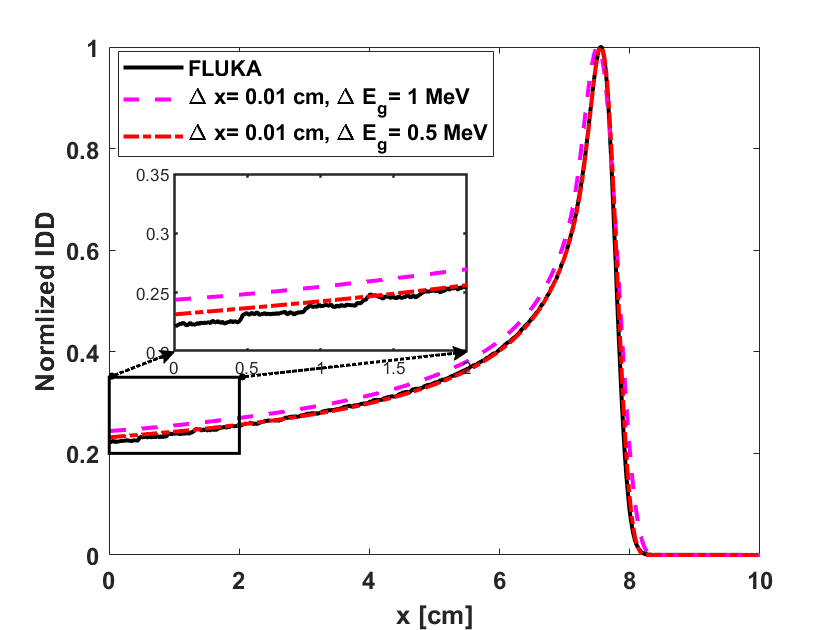}
  }
  \subfloat
  {
      \includegraphics[width=2.15in, height=2.5in]{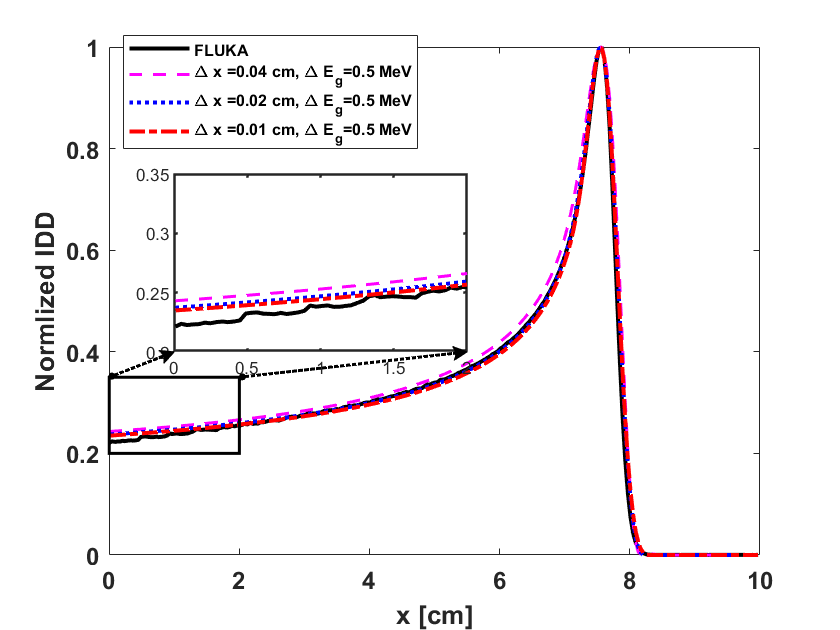}
  }
  \subfloat
  {
     \includegraphics[width=2.15in, height=2.5in]{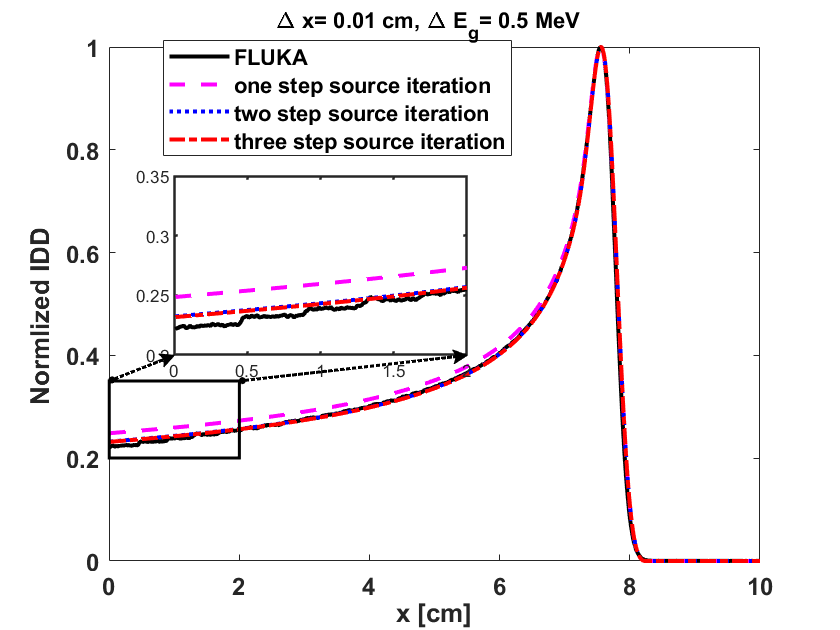}
  }	
	% \subfloat{
	% 	\includegraphics[width=0.5\textwidth]{error_1.png}
	% }
	% \subfloat{
	% 	\includegraphics[width=0.5\textwidth]{error_2.png}
	% }
	\caption{Comparison of the IDD using depth splitting deterministic method  with different spatial, energy and iteration steps.}
	\label{fig_er}
\end{figure}

This subsection aims at testing the convergence  order of our  scheme in the depth variable $x$ and energy variable $E$. The material target is water, and the energy of the incoming beam is 50 MeV. In Fig.~\ref{fig:conv}, we plot the following errors
\begin{equation} \label{error}
	\text{error}_1= \sum_{g=1}^{G}\sum_{s=1}^{N_x}\lvert  \psi_{g,\Delta x}^{1,s}-\psi_{g,\Delta x/2}^{1,s}\rvert\Delta x\Delta E_g ,\quad
	\text{error}_2= \sum_{g=1}^{G}\sum_{s=1}^{N_x}\lvert  \psi_{g,\Delta x}^{2,s}-\psi_{g,\Delta x/2}^{2,s}\rvert\Delta x\Delta E_g,
\end{equation}
for different $\Delta x=1/50\ {\rm cm}$, $1/100\ {\rm cm}$, $1/200\ {\rm cm}$, $1/400\ {\rm cm}$, $1/800\ {\rm cm}$ and different $\Delta E_g=1/2\ {\rm MeV}$, $1/4\ {\rm MeV}$, $1/8\ {\rm MeV}$, $1/16\ {\rm MeV}$, $1/32\ {\rm MeV}$. Here $\psi_{g,\Delta x}^{l,s}$ ($l=1,2$) is the discretized version of the following integration
$$
\psi_{g}^{l,s}=\int_{U,V}\int_{Y,Z}\psi^{l}(x^s,E_g)\text{d}y\text{d}z\text{d}u\text{d}v.
$$
One can observe a uniform second order accuracy for the depth and energy discretization.  
In order to demonstrate how the accuracy is effected by the spatial and energy steps, we display in Figure.~\ref{fig_er} the IDD of water with 100\ MeV incoming beam energy calculated with different spatial and energy steps. Since the clinical standard for IDD is 3 percent \cite{cli}, to achieve this accuracy, $\Delta x$=0.02\ cm, $\Delta E_g$=0.5\ MeV are required. Moreover, different iteration steps of the source iteration method are tested, one can observe that the numerical dose distribution is accurate enough to consider protons with one catastrophic scatter interaction.

\subsection{Comparison of IDD for different materials}

\begin{figure}[!h]
	\centering
	\subfloat{
		\includegraphics[width=0.5\textwidth]{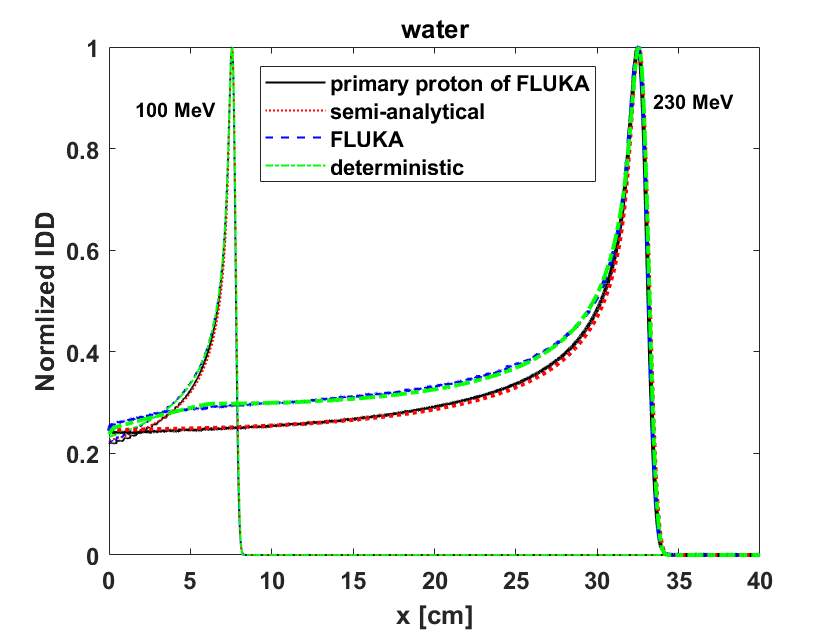}
	}
	\subfloat{
		\includegraphics[width=0.5\textwidth]{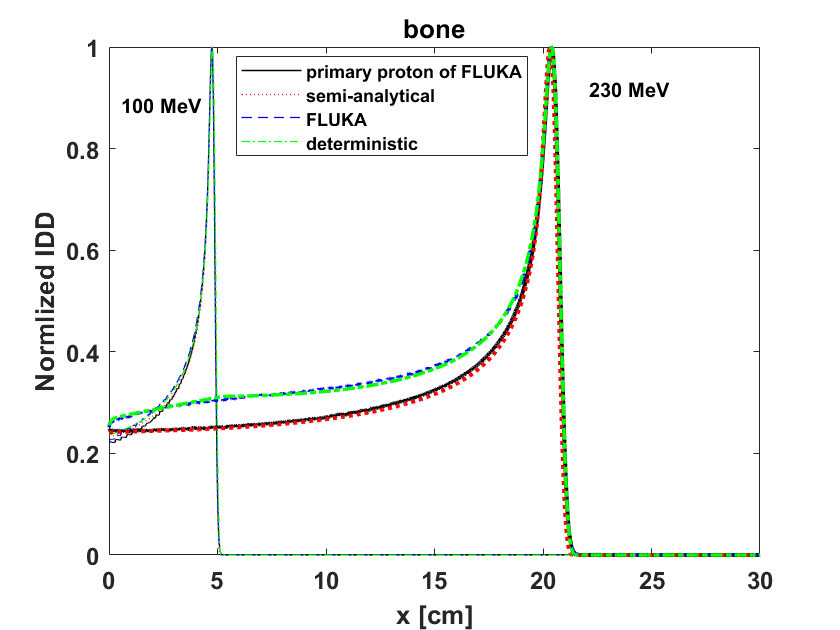}
	}
 %  \newline
 %  \subfloat{
	% 	\includegraphics[width=0.5\textwidth]{idd_ws.png}
	% }
	% \subfloat{
	% 	\includegraphics[width=0.5\textwidth]{idd_bs.png}
	% }
	\caption{Comparison of the IDD using depth splitting deterministic method  and FLUKA, and comparison of the primary IDD using semi-analytical method and FLUKA left: target material is water, right: target material is bone. Here we use 4000 $\times$ 500 mesh in $(x,E)$ for deterministic method.}
	\label{fig.1}
\end{figure}

In this subsection, we will compare the IDD  for different materials. Integrating the dose distribution laterally produces the dose-depth curve known as the integrated depth dose, i.e., 
$$
\text{IDD}(x)=\int\int_{Y,Z} D(x,y,z)dydz.
$$
    The results are shown in Fig.~\ref{fig.1}. The semi-analytical method cannot get the correct energy deposition, but its IDD matches well with the deposition from primary proton. We can see that the IDD produced by our depth splitting deterministic method almost overlap with the results given by FLUKA. In proton therapy, researchers are interested in three different positions:  the position of Bragg peak (BP), proximal-90$\%$ (P90) which is  90$\%$ dose before Bragg  peak, distal-90$\%$ (D90) which is  90$\%$ dose after Bragg peak, and distal-20$\%$ (D20). In Table.~\ref{tab:stab3}, we compare these positions in IDD  obtained by FLUKA and the depth splitting method.  We can observe that the maximum relative difference is 0.98$\%$, which shows very good agreement.

% \begin{figure}[!h]
% 	\centering
% 	\subfloat{
% 		\includegraphics[width=0.5\textwidth]{idd_wp.png}
% 	}
% 	\subfloat{
% 		\includegraphics[width=0.5\textwidth]{idd_bp.png}
% 	}
% 	\caption{Comparison of the primary proton IDD using depth splitting deterministic method, semi-analytical method  and FLUKA, left: target material is water, right: target material is bone. Here we use 4000 $\times$ 500 mesh in $(x,E)$ for deterministic method and semi-analytical method.}
% 	\label{fig.1}
% \end{figure}

\begin{table}[]
	\centering
\caption{Comparison of crucial positions in the IDD between FLUKA (F) and the deterministic method (D), here BP is the position of Bragg peak, P90 (D90) indicates the proximal-90$\%$ (distal-90$\%$) which is  90$\%$ dose before (after) Bragg  peak; D20 indicates distal-20$\%$ which is 20$\%$ dose after Bragg  peak (Unit: cm). Here we use 4000 $\times$ 500 mesh in $(x,E)$ for deterministic method.}\label{tab:stab3}	
\begin{tabular}{llllllllll}
\hline
material               & \begin{tabular}[c]{@{}l@{}}energy\\ (MeV)\end{tabular} & \begin{tabular}[c]{@{}l@{}}BP\\ F\end{tabular} & \begin{tabular}[c]{@{}l@{}}BP\\ D\end{tabular} 
& \begin{tabular}[c]{@{}l@{}}P90\\ F\end{tabular} & \begin{tabular}[c]{@{}l@{}}P90\\ D\end{tabular}&
 \begin{tabular}[c]{@{}l@{}}D90\\ F\end{tabular} & \begin{tabular}[c]{@{}l@{}}D90\\ D\end{tabular} & \begin{tabular}[c]{@{}l@{}}D20\\ F\end{tabular} & \begin{tabular}[c]{@{}l@{}}D20\\ D\end{tabular}   \\ \hline
\multirow{3}{*}{water} & 50                                                     & 2.150                                                 & 2.150                                                     & 2.083                                            & 2.075       & 2.202                                            & 2.205                              & 2.341                                               & 2.345                                  \\ \cline{2-10} 
                       & 100                                                    & 7.560                                                 & 7.560                                                         & 7.416                                            & 7.418                                            & 7.662                                            & 7.671 & 7.921                                                 & 7.945 \\ \cline{2-10} 
                       & 230                                                    & 32.460                                                & 32.540                                                     & 32.135                                           & 32.172                                        & 32.775                                           & 32.828& 33.445                                              & 33.513      \\ \hline
\multirow{3}{*}{bone}  & 50                                                     & 1.360                                                 & 1.360                                                      & 1.315                                            & 1.319                                            & 1.391                                            & 1.402& 1.473                                             & 1.488  \\ \cline{2-10} 
                       & 100                                                    & 4.740                                                 & 4.770                                                     & 4.653                                            & 4.675                                         & 4.806                                            & 4.837 & 4.971                                                & 5.005          \\ \cline{2-10} 
                       & 230                                                    & 20.360                                                & 20.390                                                 & 20.165                                           & 20.168                                         & 20.572                                           & 20.585& 21.005                                               & 21.023                \\ \hline
\end{tabular}
\end{table}

\subsection{Comparison of the spot distribution}
\begin{figure}[!htbp]    
  \centering            
  \subfloat   % 
  {
      \includegraphics[width=2.15in, height=2.5in]{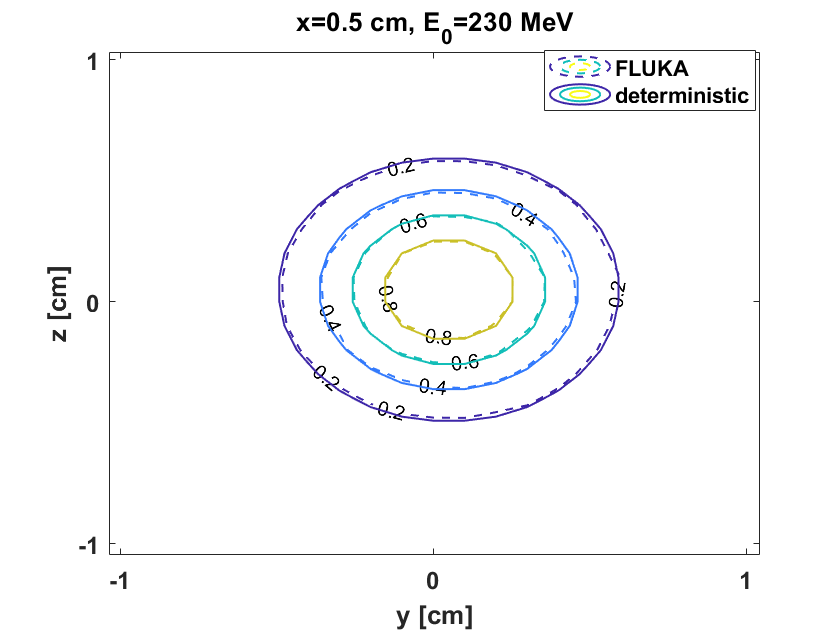}
  }
  \subfloat
  {
      \includegraphics[width=2.15in, height=2.5in]{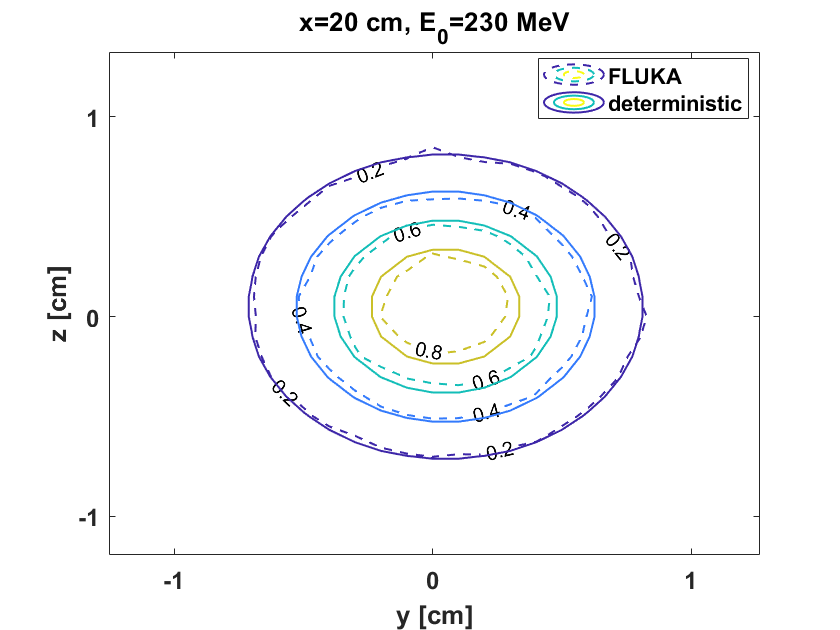}
  }
  \subfloat
  {
     \includegraphics[width=2.15in, height=2.5in]{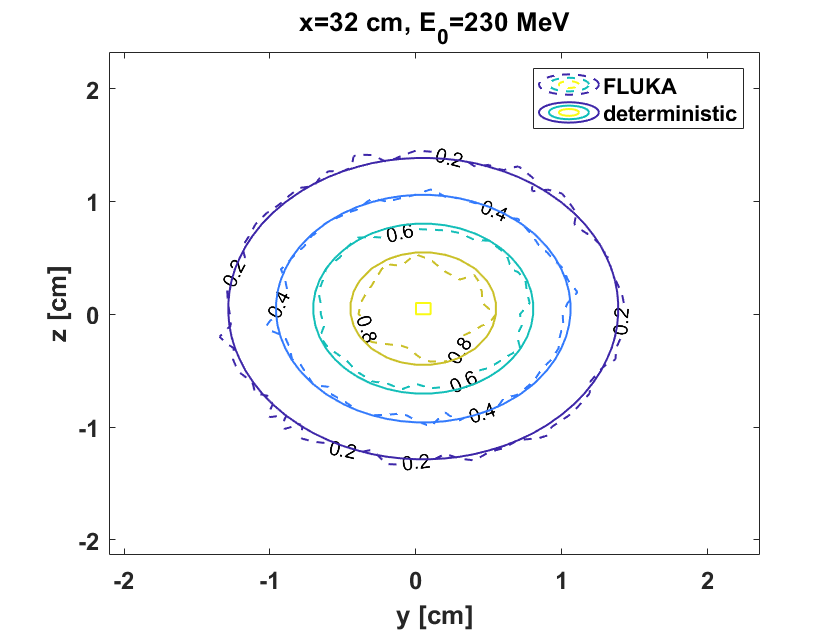}
  }
  %  \newline
  %  \subfloat   % 
  % {
  %   \includegraphics[width=2.15in, height=2.5in]{water/230/fyz1.png}
  % }
  % \subfloat
  % {
  %   \includegraphics[width=2.15in, height=2.5in]{water/230/dyz1.png}
  % }
  % \subfloat
  % {
  %     \includegraphics[width=2.15in, height=2.5in]{water/230/yzw.png}
  % }  
  % \newline
  %  \subfloat   % 
  % {
  %    \includegraphics[width=2.15in, height=2.5in]{water/230/fyz2.png}
  % }
  % \subfloat
  % {
  %     \includegraphics[width=2.15in, height=2.5in]{water/230/dyz2.png}
  % }
  % \subfloat
  % {
  %    \includegraphics[width=2.15in, height=2.5in]{water/230/yzw1.png}
  % }   
%  \newline
%   \subfloat
%  {
%  	\includegraphics[width=0.34\textwidth]{water/50/yzw5.png}
%  }
%  \subfloat
%  {
%  	\includegraphics[width=0.34\textwidth]{water/50/yzw6.png}
%  }  
%  \subfloat
%  {
%  	\includegraphics[width=0.34\textwidth]{water/50/yzw7.png}
%  }   
%   \newline
%     \subfloat
%   {
%   	\includegraphics[width=0.34\textwidth]{water/100/yzw2.png}
%   }
%   \subfloat
%   {
%   	\includegraphics[width=0.34\textwidth]{water/100/yzw3.png}
%   }  
%   \subfloat
%   {
%   	\includegraphics[width=0.34\textwidth]{water/100/yzw4.png}
%   }   
  \caption{Comparison of the spot distributions using the deterministic method and FLUKA in water.}    % 
  \label{fig.2}            % 
\end{figure}

\begin{figure}[!htbp]    
	\centering            
	\subfloat   % 
	{
		\includegraphics[width=2.15in, height=2.5in]{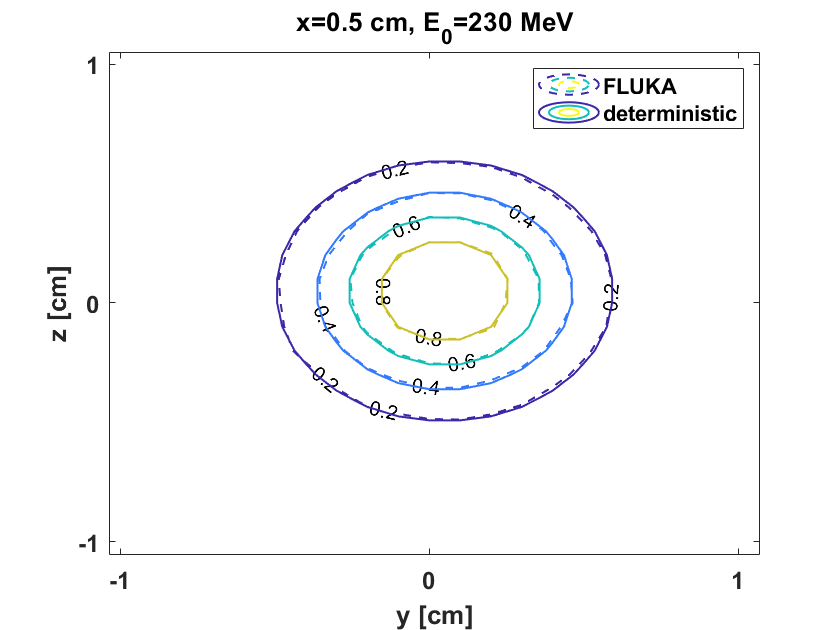}
	}
	\subfloat
	{
		\includegraphics[width=2.15in, height=2.5in]{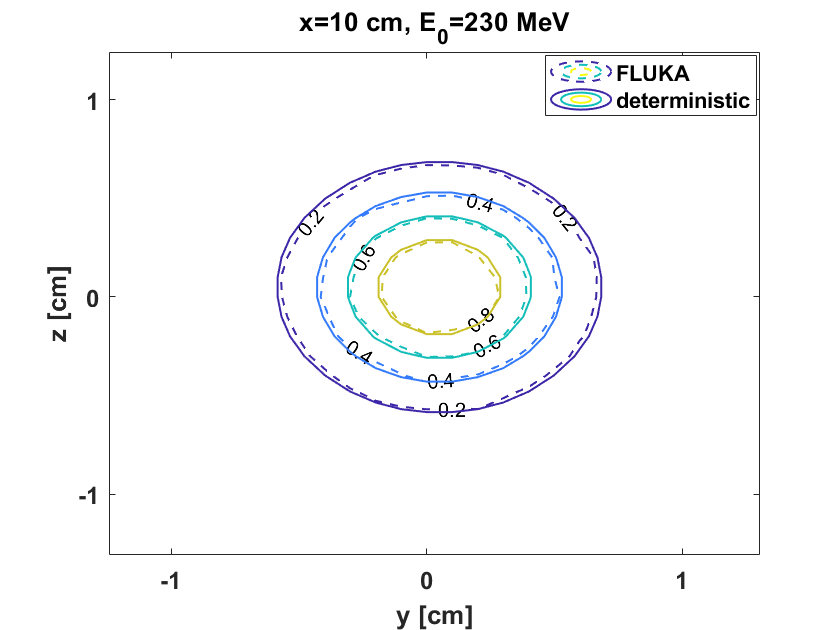}
	}
	\subfloat
	{
		\includegraphics[width=2.15in, height=2.5in]{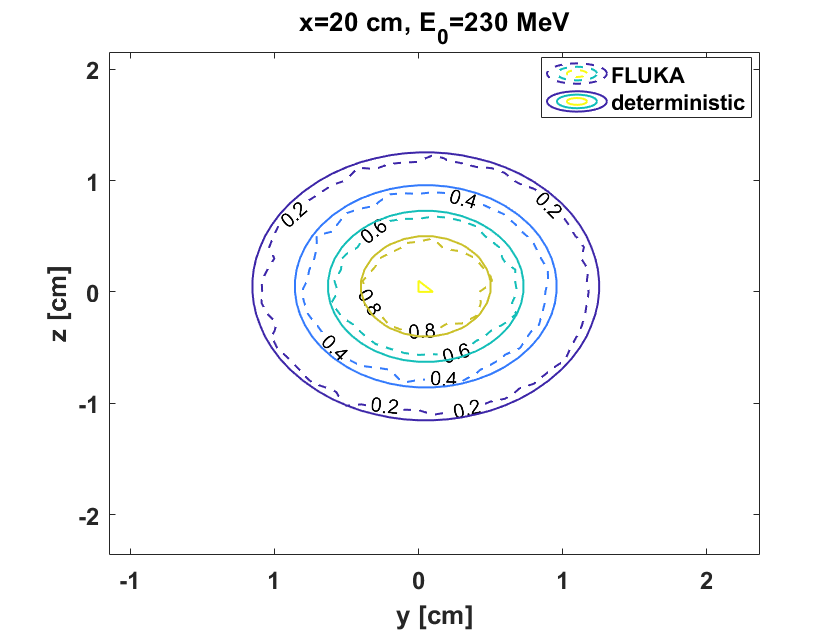}
	}
	% \newline
	% \subfloat   % 
	% {
	% 	\includegraphics[width=2.15in, height=2.5in]{bone/230/fxy2.png}
	% }
	% \subfloat
	% {
	% 	\includegraphics[width=2.15in, height=2.5in]{bone/230/dxy2.png}
	% }
	% \subfloat
	% {
	% 	\includegraphics[width=2.15in, height=2.5in]{bone/230/yzb8.png}
	% }  
	% \newline
	% \subfloat   % 
	% {
	% 	\includegraphics[width=2.15in, height=2.5in]{bone/230/fxy3.png}
	% }
	% \subfloat
	% {
	% 	\includegraphics[width=2.15in, height=2.5in]{bone/230/dxy3.png}
	% }
	% \subfloat
	% {
	% 	\includegraphics[width=2.15in, height=2.5in]{bone/230/yzb9.png}
	% }   
%	\newline
%	\subfloat
%	{
%		\includegraphics[width=0.34\textwidth]{bone/50/yzw1.png}
%	}
%	\subfloat
%	{
%		\includegraphics[width=0.34\textwidth]{bone/50/yzw2.png}
%	}  
%	\subfloat
%	{
%		\includegraphics[width=0.34\textwidth]{bone/50/yzw3.png}
%	}   
%	\newline
%	\subfloat
%	{
%		\includegraphics[width=0.34\textwidth]{bone/100/yzb4.png}
%	}
%	\subfloat
%	{
%		\includegraphics[width=0.34\textwidth]{bone/100/yzb5.png}
%	}  
%	\subfloat
%	{
%		\includegraphics[width=0.34\textwidth]{bone/100/yzb6.png}
%	}   
	\caption{Comparison of the spot distributions using the deterministic method and FLUKA in bone.}    % 
	\label{fig.3}            % 
\end{figure}

In this subsection, we will compare the spot distribution for different materials. The spot distribution is the dose distribution on the $YZ$ plane at a given depth. Here, we have selected three different positions for different $E_0$: at the entrance, in the middle, and near the peak. We show the results of $E_0=230$ MeV. The numerical results are presented in Figure.~\ref{fig.2} with the target material being water. The numerical results are gathered in Figure.~\ref{fig.3} with the target material being bone. Moreover, we can calculate the standard deviation $\sigma$ of the lateral distribution, and the results for three different materials (as in Table 1) are listed in Table.~\ref{tab:stab9_w}--\ref{tab:stab9_a}. Here, we can see that our results show very good agreement with the reference solution. 
%\textcolor{blue}{Moreover, we can also see that the  accuracy of the semi-analytical method in lower than our method.}
% \usepackage{multirow}
\begin{table}[]
\centering
\caption{Comparison of $\sigma$ of the lateral distributions at three different depths in water for beams with different injected energy levels. EN is at the entrance $x=0$, MI is in the middle, and NP is near the peak. Here F is for the results obtained by FLUKA and D indicates results obtained by the deterministic method(Unit: cm).}\label{tab:stab9_w}
\begin{tabular}{llllllllll}
\hline
Energy (MeV)           & \multicolumn{3}{l}{50}   & \multicolumn{3}{l}{100}  & \multicolumn{3}{l}{230}  \\
\hline
\multirow{2}{*}{Depth} & EN     & MI     & NP     & EN     & MI     & NP     & EN     & MI     & NP     \\
                       & 0.5    & 1.5    & 2.1    & 0.5    & 3.5    & 7.5    & 0.5    & 20     & 32     \\
\hline
F                      & 0.2971 & 0.3002 & 0.3035 & 0.2977 & 0.3006 & 0.3502 & 0.2967 & 0.4282 & 0.7541 \\
\hline
D                      & 0.3000 & 0.3018 & 0.3049 & 0.3000 & 0.3049 & 0.3543 & 0.3000 & 0.4306 & 0.7570\\
\hline
%\textcolor{blue}{S}                      & \textcolor{blue}{0.3001} & \textcolor{blue}{0.3016} &\textcolor{blue}{ 0.3047} & \textcolor{blue}{0.3000 }& \textcolor{blue}{0.3034} & \textcolor{blue}{0.3355} & \textcolor{blue}{0.3000} & \textcolor{blue}{0.4121} & \textcolor{blue}{0.7409}\\
%\hline
\end{tabular}
\end{table}
\begin{table}[]
\centering
\caption{Comparison of $\sigma$ of the lateral distributions at three different depths in bone for beams with different injected energy levels. EN is at the entrance $x=0$, MI is in the middle, and NP is near the peak. Here F is for the results obtained by FLUKA and D indicates results obtained by the deterministic method(Unit: cm).}\label{tab:stab9_b}
\begin{tabular}{llllllllll}
\hline
Energy (MeV)           & \multicolumn{3}{l}{50}   & \multicolumn{3}{l}{100}  & \multicolumn{3}{l}{230}  \\
\hline
\multirow{2}{*}{Depth} & EN     & MI     & NP     & EN     & MI     & NP     & EN     & MI     & NP     \\
                       & 0.5    & 1.0    & 1.3    & 0.5    & 2.0    & 4.7    & 0.5    & 10.0   & 20.0   \\
\hline
F                      & 0.2976 & 0.2991 & 0.2996 & 0.2964 & 0.2979 & 0.3325 & 0.2998 & 0.3483 & 0.6225 \\
\hline
D                      & 0.3000 & 0.3014 & 0.3030 & 0.3000 & 0.3021 & 0.3341 & 0.3000 & 0.3491 & 0.6277\\
\hline
%\textcolor{blue}{S}                      & \textcolor{blue}{0.3001} & \textcolor{blue}{0.3012} &\textcolor{blue}{ 0.3031} & \textcolor{blue}{0.3000 }& \textcolor{blue}{0.3011} & \textcolor{blue}{0.3284} & \textcolor{blue}{0.3000} & \textcolor{blue}{0.3374} & \textcolor{blue}{0.6012}\\
%\hline
\end{tabular}
\end{table}
\begin{table}[]
\centering
\caption{Comparison of $\sigma$ of the lateral distributions at three different depths in air for beams with different injected energy levels. EN is at the entrance $x=0$, MI is in the middle, and NP is near the peak. Here F is for the results obtained by FLUKA and D indicates results obtained by the deterministic method(Unit: cm).}\label{tab:stab9_a}
\begin{tabular}{llllllllll}
\hline
Energy (MeV)           & \multicolumn{3}{l}{50}   & \multicolumn{3}{l}{100}  & \multicolumn{3}{l}{230}  \\
\hline
\multirow{2}{*}{Depth} & $x_1$     & $x_2$     & $x_3$     & $x_1$     & $x_2$     & $x_3$     & $x_1$     & $x_2$     & $x_3$       \\
                       & 0.5    & 20   & 40   & 0.5    & 20    & 40    & 0.5    & 20     & 40     \\
\hline
F                      & 0.2972 & 0.2981 & 0.3070 & 0.2957 & 0.2985 & 0.3012 & 0.2964 & 0.2977 & 0.2983 \\
\hline
D                      & 0.3000 & 0.3012 & 0.3098 & 0.3000 & 0.3003 & 0.3026 & 0.3000 & 0.3000 & 0.3005\\
\hline
\end{tabular}
\end{table}
\subsection{Comparison of longitudinal data}
\begin{figure}[!htbp]    
  \centering            
  \subfloat   % 
  {
      \includegraphics[width=2.15in, height=2.5in]{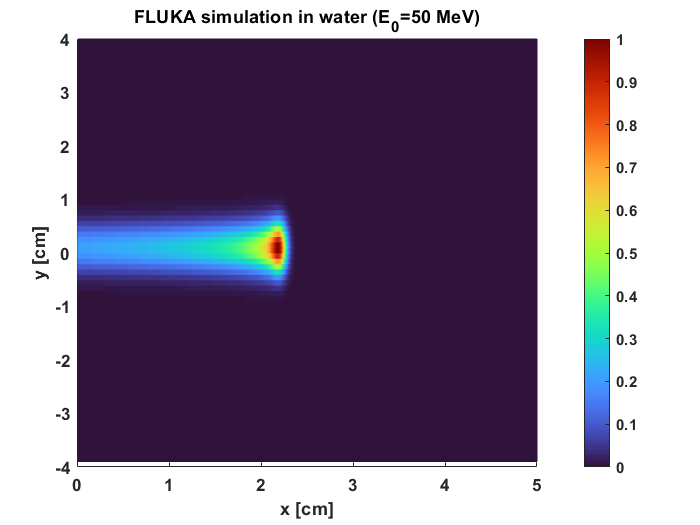}
  }
  \subfloat
  {
      \includegraphics[width=2.15in, height=2.5in]{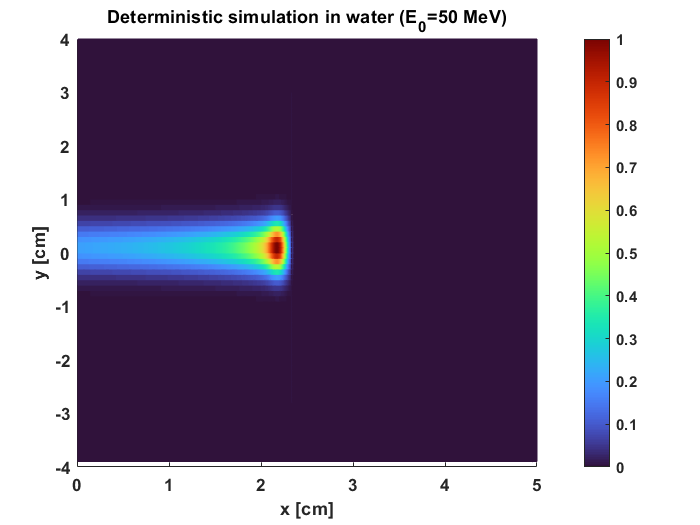}
  }
   \subfloat
  {
  	\includegraphics[width=2.15in, height=2.5in]{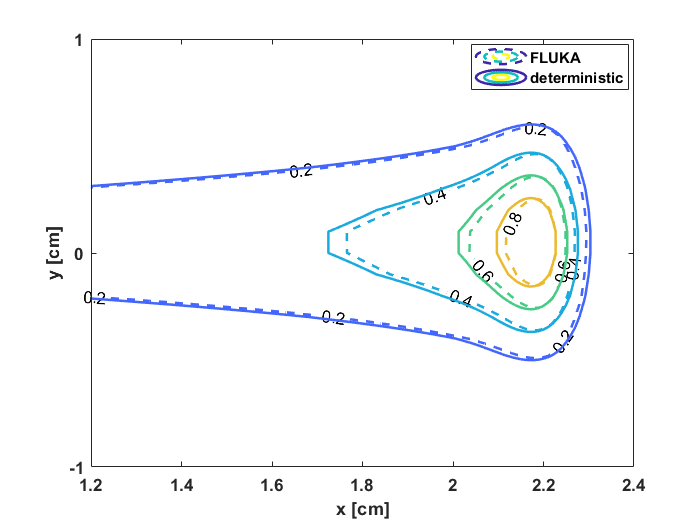}
  }
   \newline
   \subfloat   % 
  {
    \includegraphics[width=2.15in, height=2.5in]{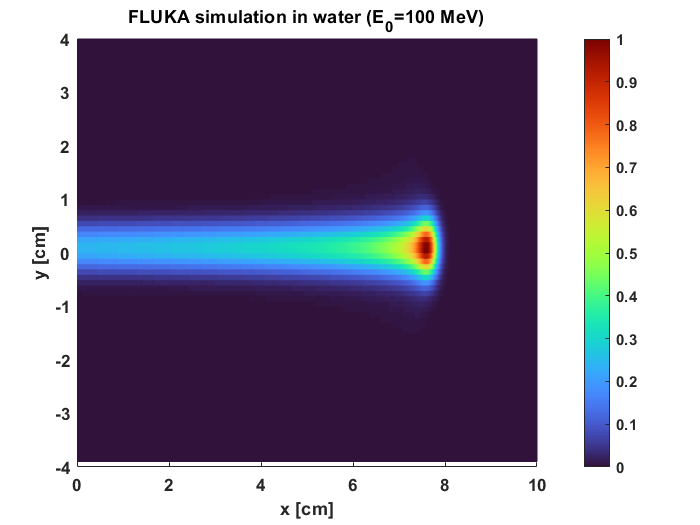}
  }
  \subfloat
  {
    \includegraphics[width=2.15in, height=2.5in]{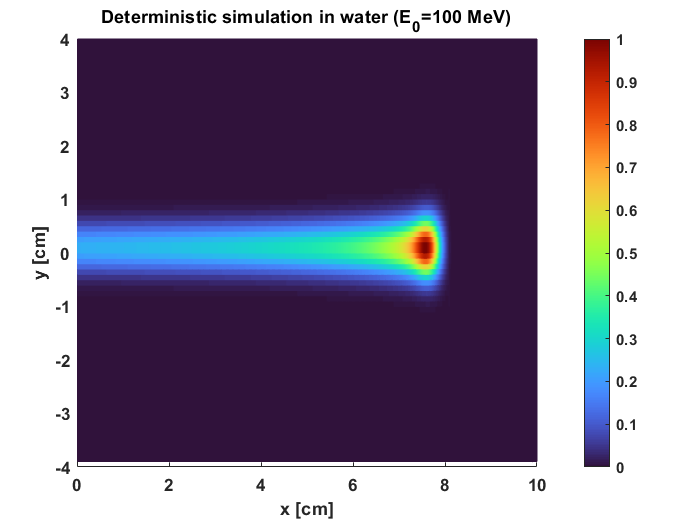}
  }
  \subfloat
  {
  	\includegraphics[width=2.15in, height=2.5in]{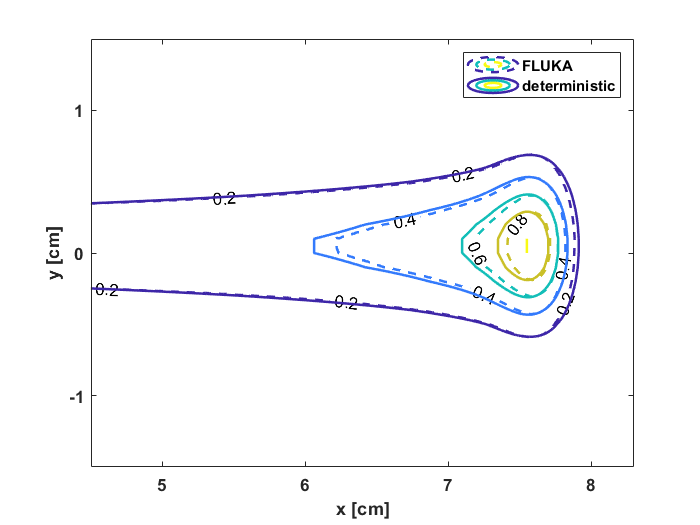}
  }
  \newline
   \subfloat   % 
  {
     \includegraphics[width=2.15in, height=2.5in]{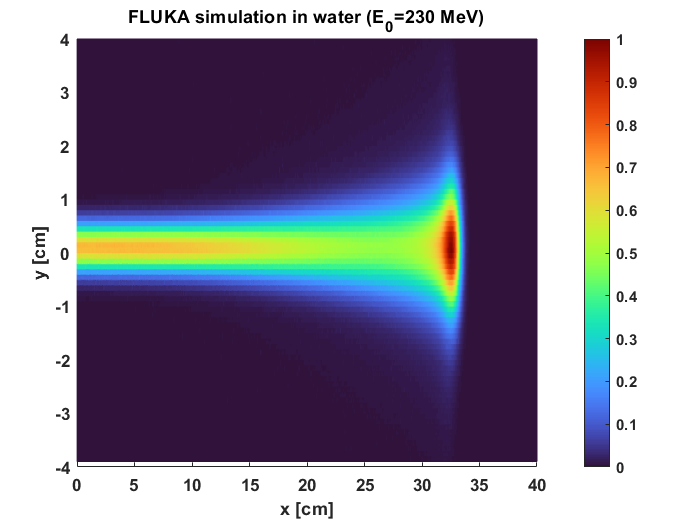}
  }
  \subfloat
  {
      \includegraphics[width=2.15in, height=2.5in]{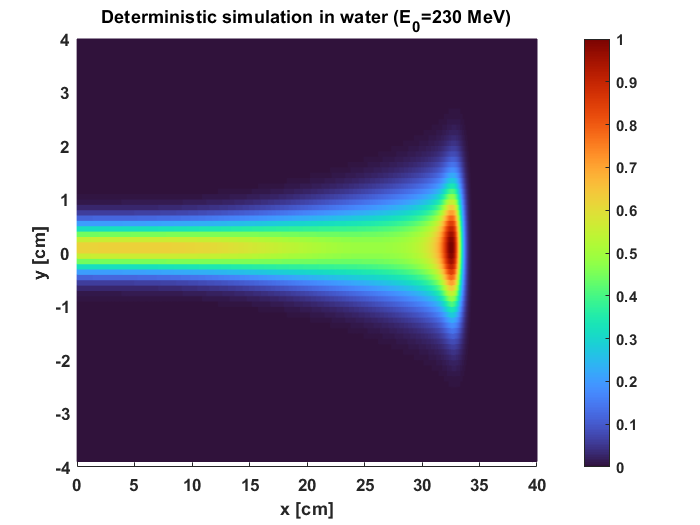}
  }
  \subfloat
  {
  	\includegraphics[width=2.15in, height=2.5in]{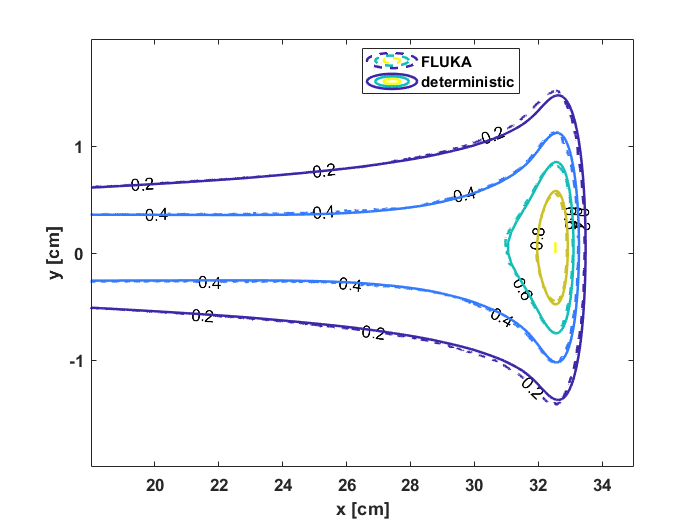}
  }
  \caption{Comparison of the longitudinal data using the deterministic method  and FLUKA for water with different energy.}    % 
  \label{fig.4}            % 
\end{figure}

\begin{figure}[!htbp]    
  \centering            
  \subfloat   % 
  {
      \includegraphics[width=2.15in, height=2.5in]{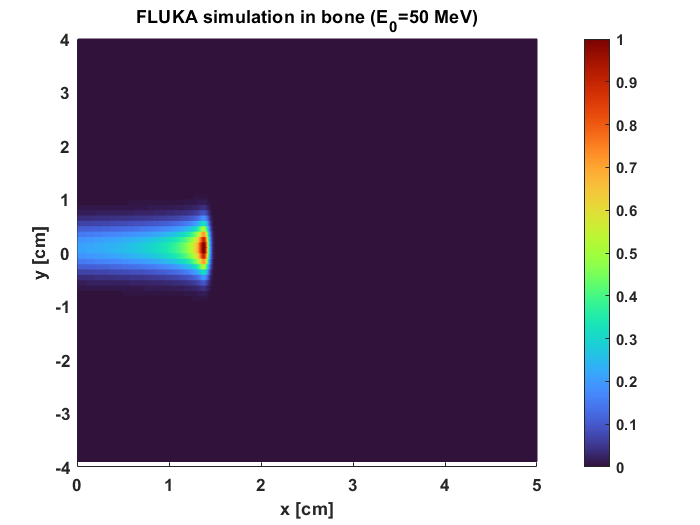}
  }
  \subfloat
  {
      \includegraphics[width=2.15in, height=2.5in]{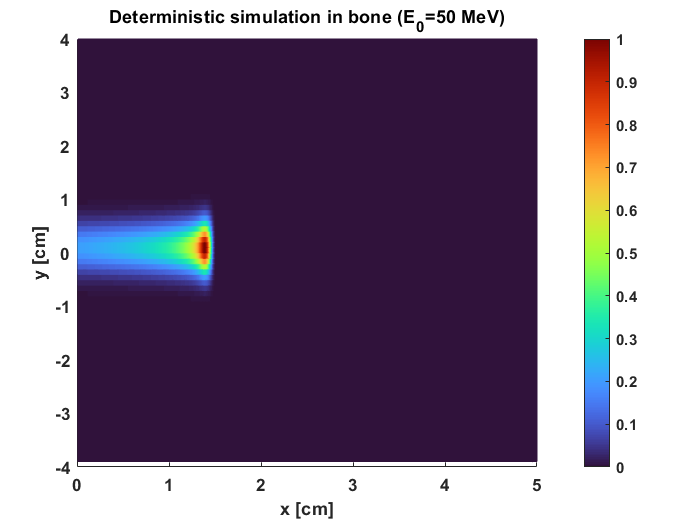}
  }
   \subfloat
  {
  	\includegraphics[width=2.15in, height=2.5in]{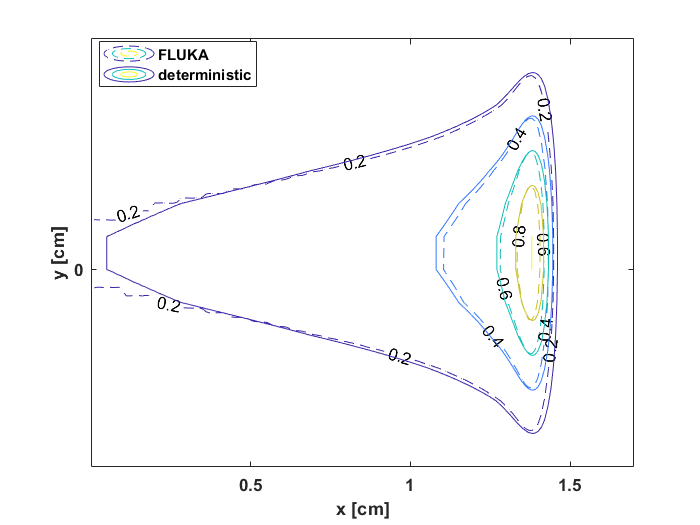}
  }
   \newline
   \subfloat   % 
  {
    \includegraphics[width=2.15in, height=2.5in]{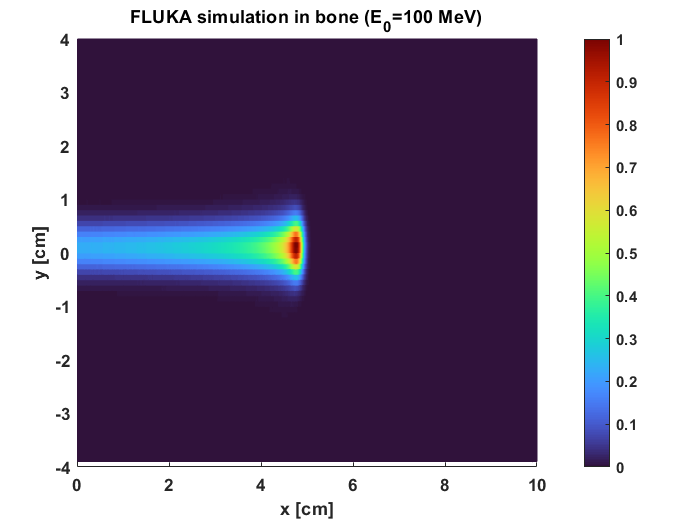}
  }
  \subfloat
  {
    \includegraphics[width=2.15in, height=2.5in]{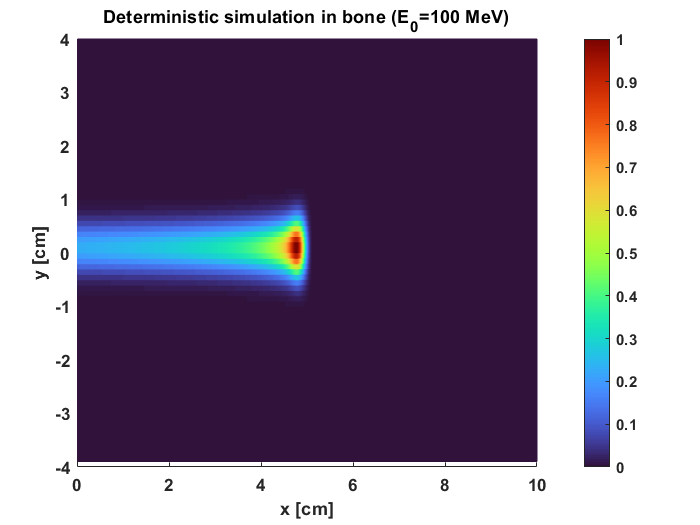}
  }
  \subfloat
  {
  	\includegraphics[width=2.15in, height=2.5in]{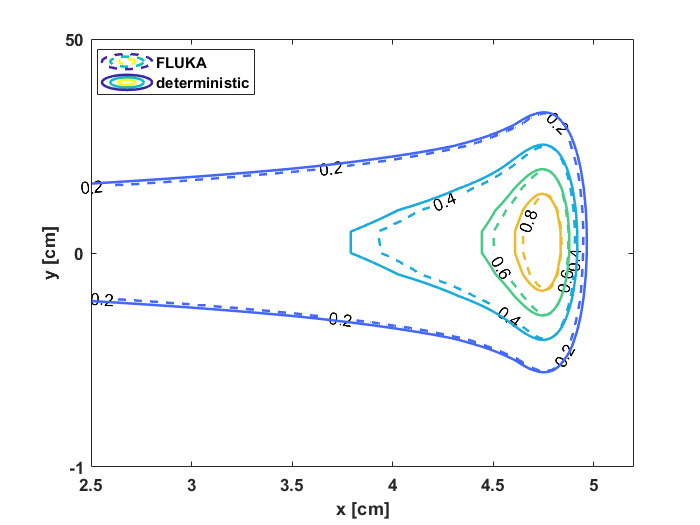}
  }
  \newline
   \subfloat   % 
  {
     \includegraphics[width=2.15in, height=2.5in]{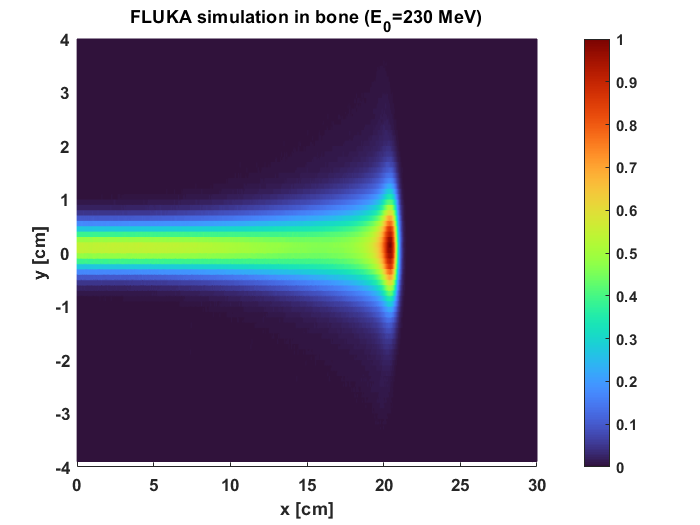}
  }
  \subfloat
  {
      \includegraphics[width=2.15in, height=2.5in]{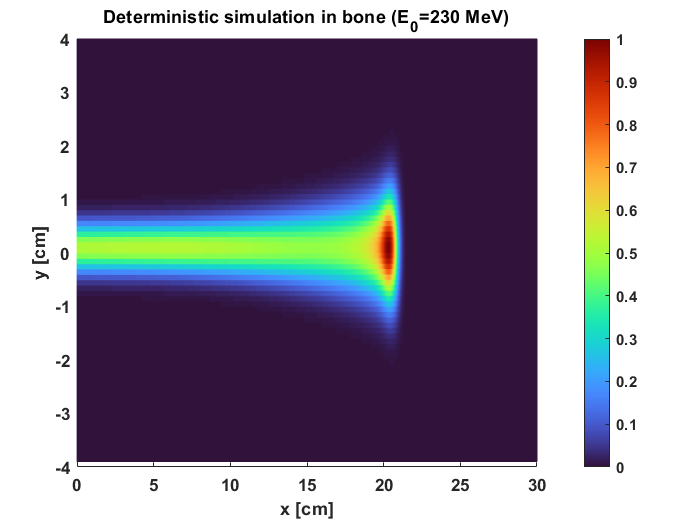}
  }
  \subfloat
  {
  	\includegraphics[width=2.15in, height=2.5in]{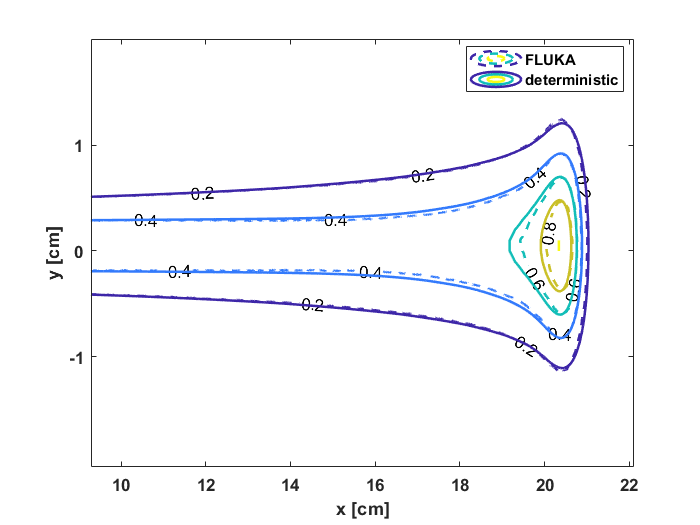}
  }
  \caption{Comparison of the longitudinal data using the deterministic method  and FLUKA for bone with different energy.}    % 
  \label{fig.5}            % 
\end{figure}

In this subsection, we will compare the longitudinal  data (LD)  for different materials. The longitudinal  data can be obtained by integrating the dose distribution in $z$-direction, i.e., 
$$
\text{LD}(x,y)=\int_{Z} D(x,y,z)dz.
$$
The numerical results are presented in Figure.~\ref{fig.4} with the target material being water. The numerical results are gathered in Figure.~\ref{fig.5} with the target material being bone. Moreover, we compare the contour lines between
the deterministic method and FLUKA. We can see that our results show  good agreement with the reference solution.

% \begin{figure}[!htbp]    
%   \centering            
%   \subfloat   % 
%   {
%       \includegraphics[width=2.15in, height=2.5in]{fp/water/e_1.png}
%   }
%   \subfloat
%   {
%       \includegraphics[width=2.15in, height=2.5in]{fp/water/e_2.png}
%   }
%    \subfloat
%   {
%   	\includegraphics[width=2.15in, height=2.5in]{fp/water/e_3.png}
%   }
%   \caption{error.}    %           % 
% \end{figure}
\subsection{A test case with non-homogeneous materials}
\begin{figure}[ht]
	\centering
	\includegraphics[width=0.6\textwidth]{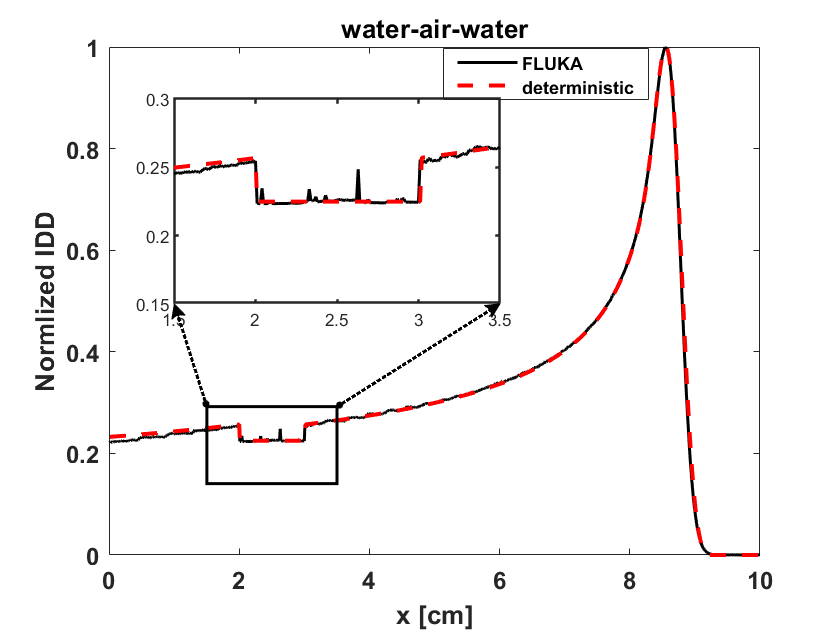}
	\caption{Comparison of the IDD using splitting deterministic method  and FLUKA, and the material water and air are co-existed.}
	\label{fig:larsen}
\end{figure}

\begin{figure}[!htbp]    
  \centering            
  \subfloat   % 
  {
      \includegraphics[width=2.15in, height=2.5in]{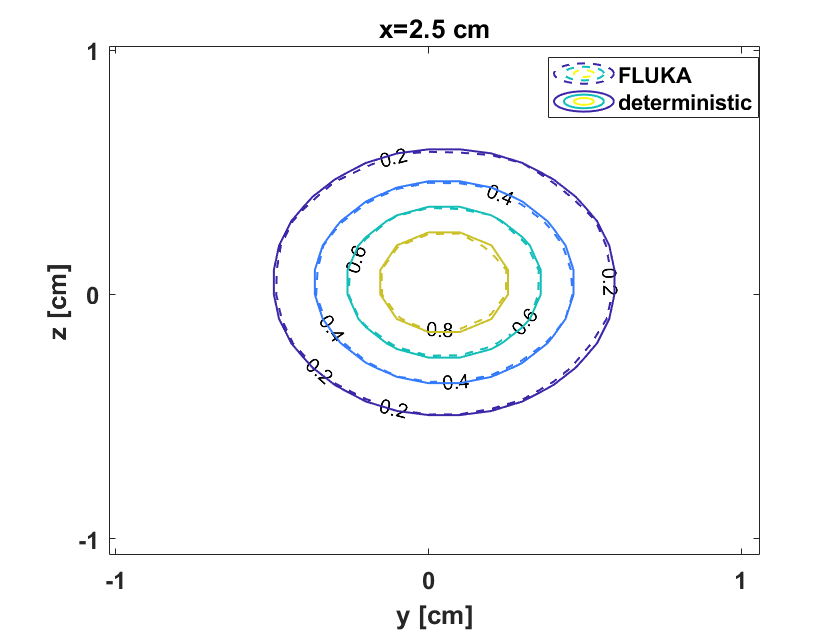}
  }
  \subfloat
  {
      \includegraphics[width=2.15in, height=2.5in]{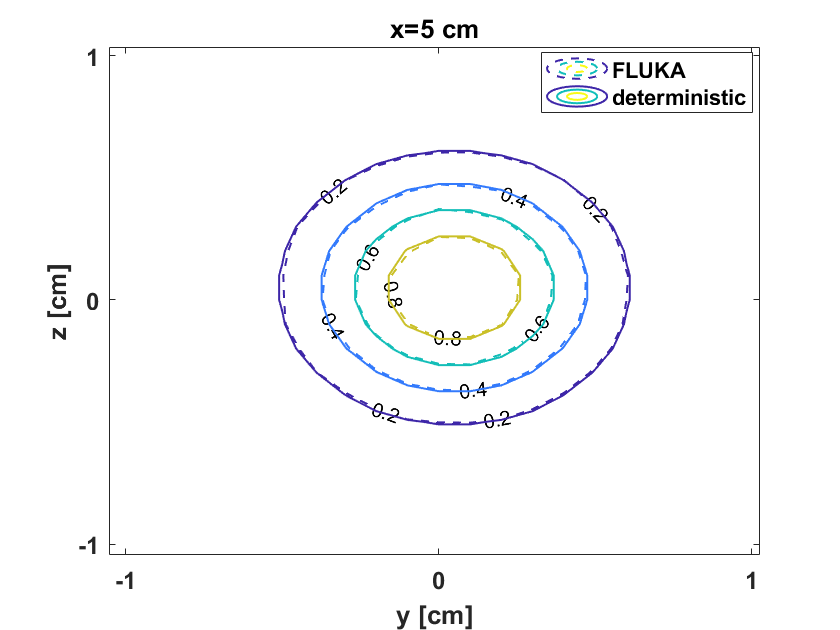}
  }
   \subfloat
  {
  	\includegraphics[width=2.15in, height=2.5in]{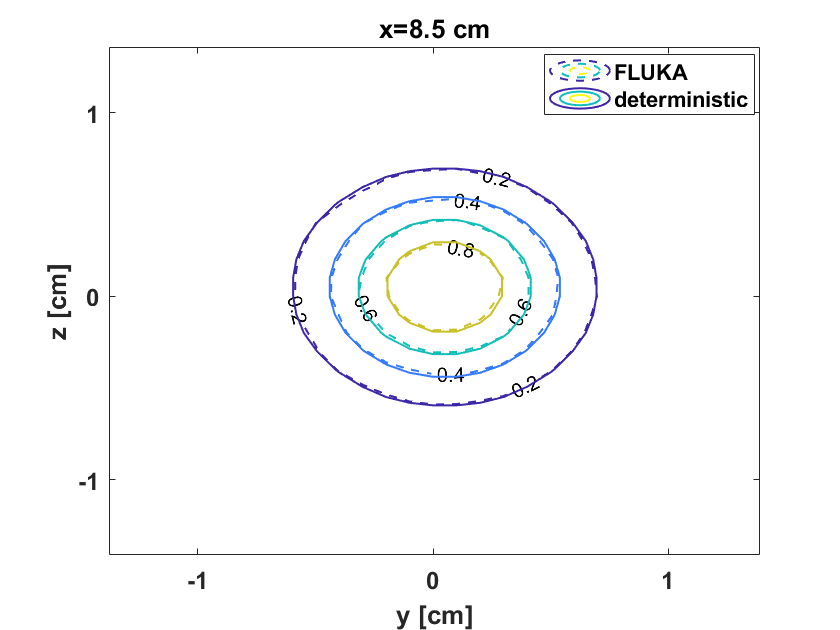}
  }
  %  \newline
  %  \subfloat   % 
  % {
  %   \includegraphics[width=2.15in, height=2.5in]{fyz25/fyz5.png}
  % }
  % \subfloat
  % {
  %   \includegraphics[width=2.15in, height=2.5in]{fyz25/dyz5.png}
  % }
  % \subfloat
  % {
  % 	\includegraphics[width=2.15in, height=2.5in]{fyz25/yz5.png}
  % }
  % \newline
  %  \subfloat   % 
  % {
  %    \includegraphics[width=2.15in, height=2.5in]{fyz25/fyz85.png}
  % }
  % \subfloat
  % {
  %     \includegraphics[width=2.15in, height=2.5in]{fyz25/dyz85.png}
  % }
  % \subfloat
  % {
  % 	\includegraphics[width=2.15in, height=2.5in]{fyz25/yz85.png}
  % }
  \caption{Comparison of the lateral data using the deterministic method  and FLUKA for water-air-water.}    % 
  \label{fig.6}            % 
\end{figure}

\begin{figure}[!htbp]    
  \centering            
  \subfloat   % 
  {
      \includegraphics[width=2.15in, height=2.5in]{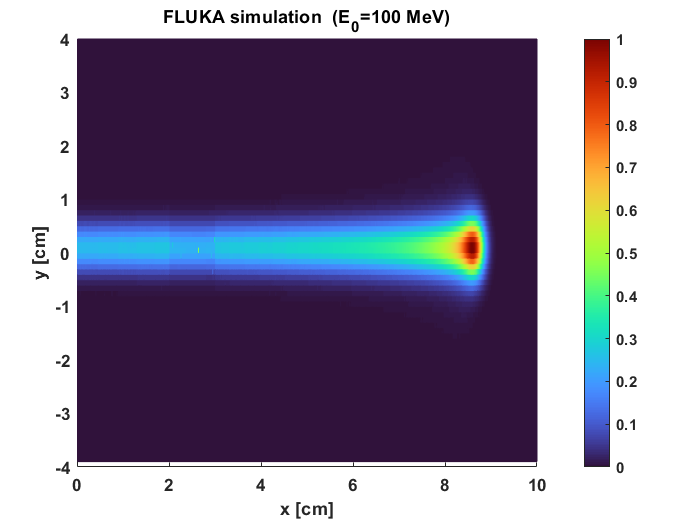}
  }
  \subfloat
  {
      \includegraphics[width=2.15in, height=2.5in]{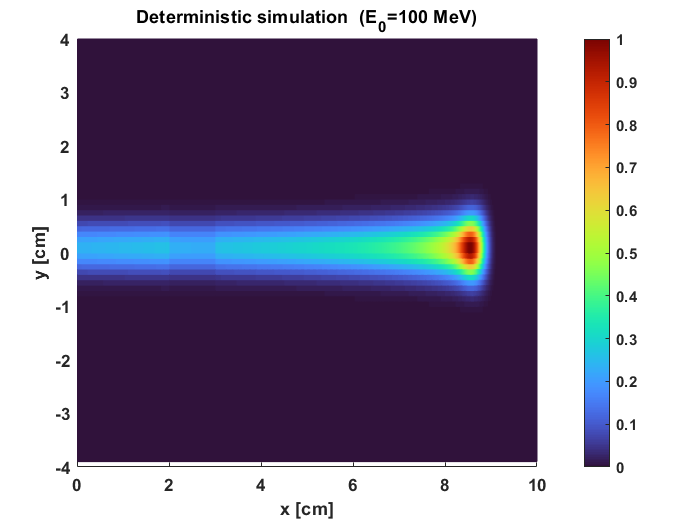}
  }
   \subfloat
  {
  	\includegraphics[width=2.15in, height=2.5in]{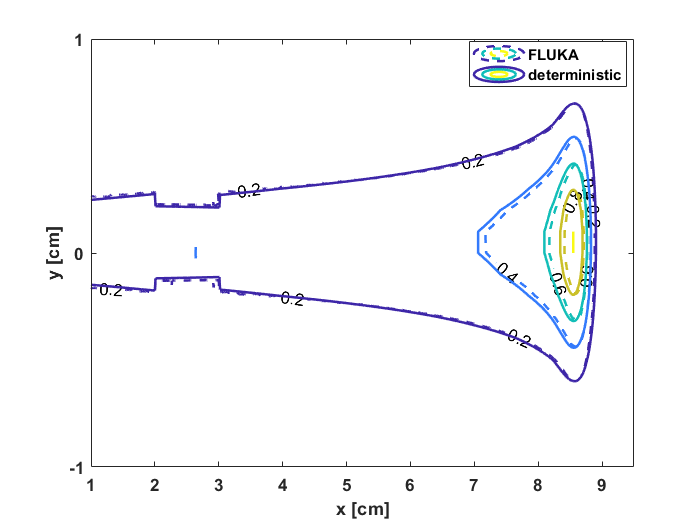}
  }
  \caption{Comparison of the longitudinal data using the deterministic method  and FLUKA for water-air-water.}    % 
  \label{fig.7}            % 
\end{figure}

In this section, we test a case with non-homogeneous medium. The main material in this case is water, a 1 cm layer of  air is placed at $x$ = 2 cm, and the initial energy $E_0$=100 MeV. The numerical results including IDD, spot distribution and longitudinal data are gathered in Fig.~\ref{fig:larsen}--Fig.~\ref{fig.7}. We can observe unrealistic oscillations inside air in the IDD obtained by FLUKA. These oscillations are due to the statistical noises of MC method, which demonstrate the advantage of using deterministic method. In other places, a good agreement between the deterministic splitting method and FLUKA can be observed.

\subsection{A test case with heterogeneous materials}
\begin{figure}[ht]
	\centering
	\includegraphics[width=0.8\textwidth]{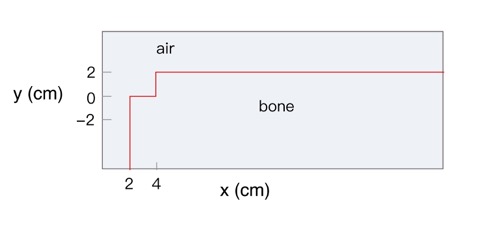}
	\caption{The top view of the  target tank.}
	\label{fig:top1}
\end{figure}

\begin{figure}[ht]
	\centering
	\includegraphics[width=0.6\textwidth]{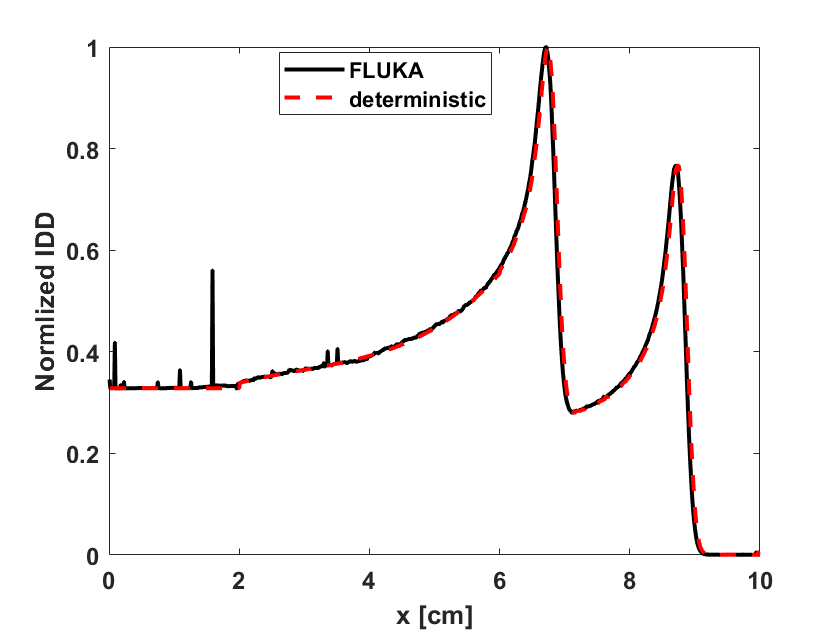}
	\caption{Comparison of the IDD using the deterministic method  and FLUKA for the case with heterogeneous medium.}
	\label{fig:h1}
\end{figure}

\begin{figure}[!htbp]    
  \centering            
  \subfloat   % 
  {
      \includegraphics[width=2.15in, height=2.5in]{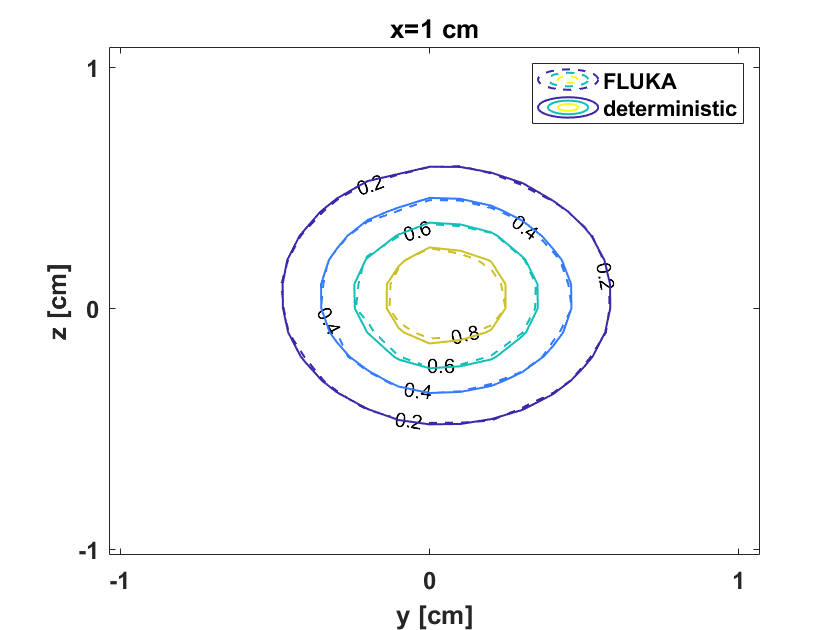}
  }
  \subfloat
  {
      \includegraphics[width=2.15in, height=2.5in]{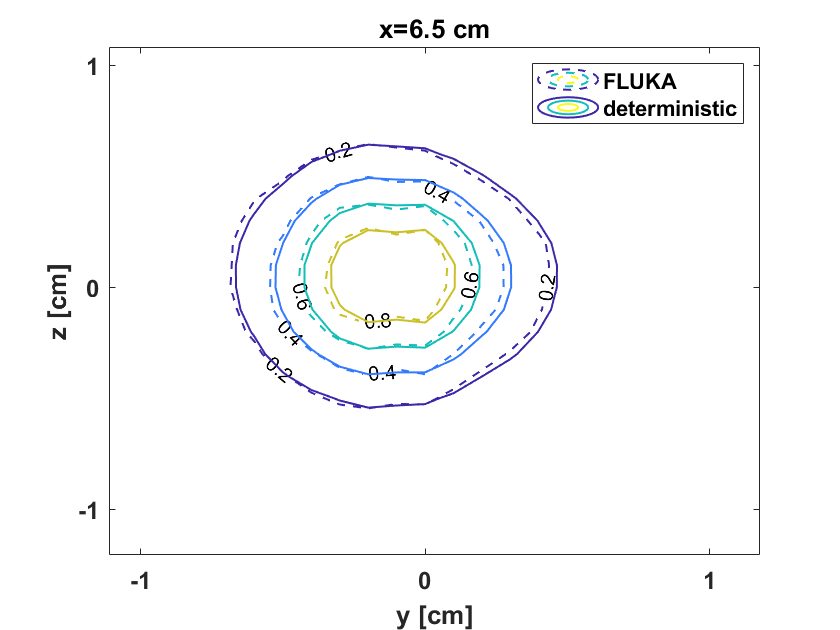}
  }
   \subfloat
  {
  	\includegraphics[width=2.15in, height=2.5in]{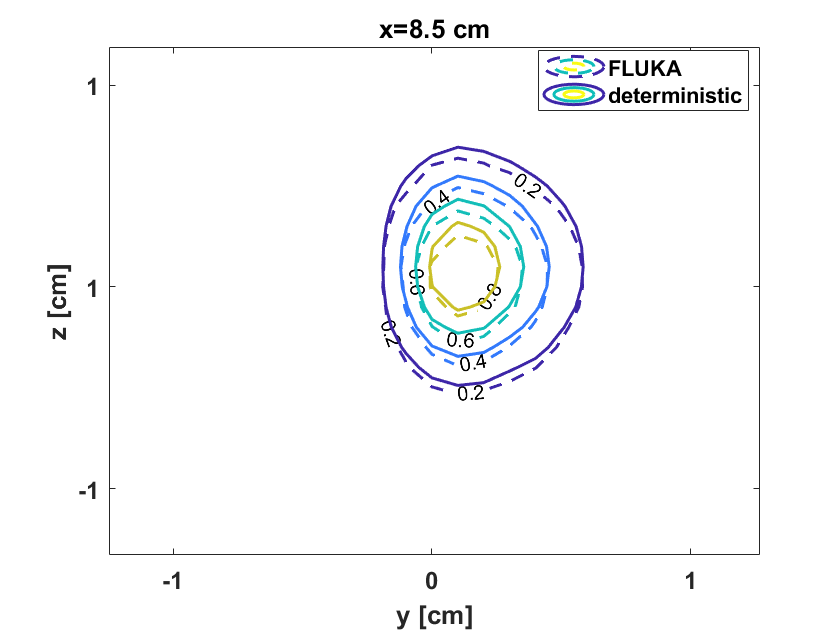}
  }
  %  \newline
  %  \subfloat   % 
  % {
  %   \includegraphics[width=2.15in, height=2.5in]{hete/fyz46.png}
  % }
  % \subfloat
  % {
  %   \includegraphics[width=2.15in, height=2.5in]{hete/dyz46.png}
  % }
  % \subfloat
  % {
  % 	\includegraphics[width=2.15in, height=2.5in]{hete/y46.png}
  % }
  % \newline
  %  \subfloat   % 
  % {
  %    \includegraphics[width=2.15in, height=2.5in]{hete/fyz75.png}
  % }
  % \subfloat
  % {
  %     \includegraphics[width=2.15in, height=2.5in]{hete/dyz75.png}
  % }
  % \subfloat
  % {
  % 	\includegraphics[width=2.15in, height=2.5in]{hete/y75.png}
  % }
  \caption{Comparison of the lateral data using the deterministic method  and FLUKA for heterogeneous case.}    % 
  \label{fig:h2}            % 
\end{figure}

\begin{figure}[!htbp]    
  \centering            
  \subfloat   % 
  {
      \includegraphics[width=2.15in, height=2.5in]{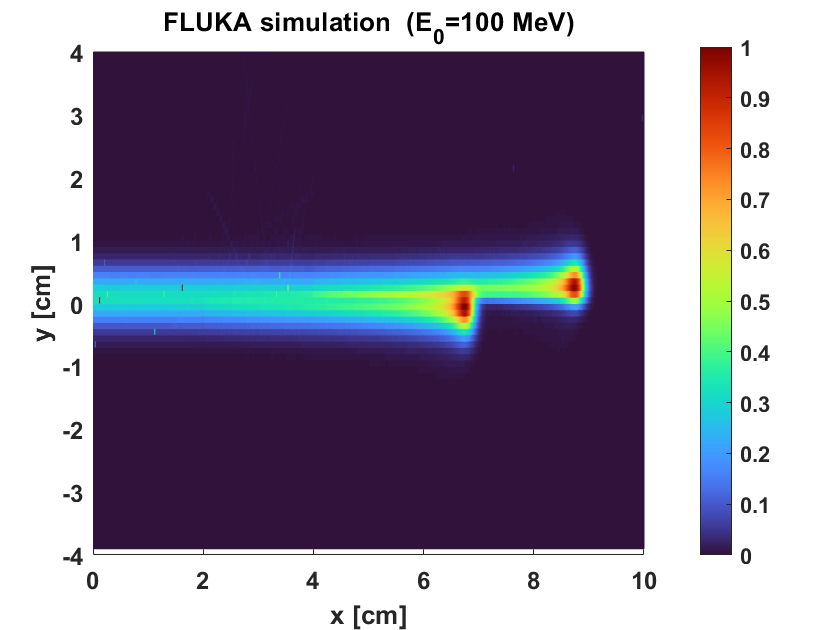}
  }
  \subfloat
  {
      \includegraphics[width=2.15in, height=2.5in]{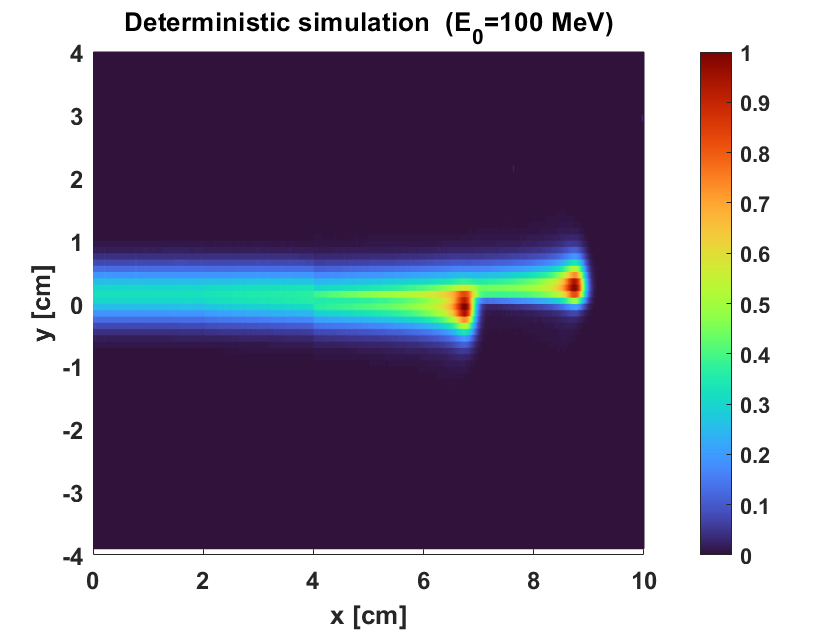}
  }
   \subfloat
  {
  	\includegraphics[width=2.15in, height=2.5in]{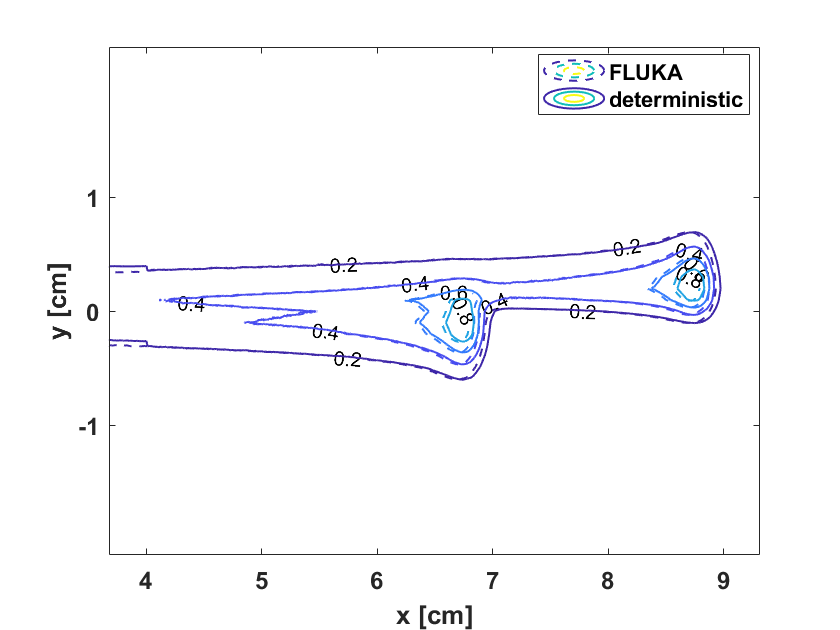}
  }
  \caption{Comparison of the longitudinal data using the deterministic method  and FLUKA for heterogeneous case.}    % 
  \label{fig:h3}            % 
\end{figure}

Finally, we test a case with heterogeneous medium with no layer structure. The top view of the  tank is shown in Fig.~\ref{fig:top1}, and the initial energy $E_0$=100 MeV. The numerical results  including IDD, spot distribution and longitudinal data are gathered in Fig.~\ref{fig:h1}--Fig.~\ref{fig:h3}. In the results, a good agreement between the deterministic splitting method and FLUKA can be observed. Moreover, the  lateral heterogeneities cannot be solved by the the algorithm proposed in \cite{burlacu2023deterministic} even for primary proton, and our algorithm can account for lateral heterogeneities case.


\section{Conclusion and discussions}
The linear boltzmann model for the proton distribution is a 6-dimensional PDE. A deterministic depth splitting method which is easy to parallelize is proposed in this paper. Firstly, second order DG method is used for the energy discretization, in which the solution is approximated by piecewise linear functions inside each energy group. The  coefficients of the linear functions inside each group satisfy a 5-dimensional equation. Then second order strang splitting method is employed for the variable representing depth. Finite difference method is used for the velocity discretization. Since the meshes for velocity are equally spaced, Fourier method can be applied to update the diffusion operator in velocity. Second order MUSCL scheme is used to discretize the linear transport in spatial variable.

In this paper, we consider not only inelastic scattering and elastic scattering, but also catastrophic scattering. We numerically verified that by using the two-step source iteration method, i.e. by considering only the primary protons and secondary protons, the obtained energy deposition matches well with the full MC method. Moreover, for mono-directional proton beams whose energy and proton movement directions are concentrated at the nozzle, beam will not expand too much in the transverse direction and the Fermi-Eyges approximation for elastic scattering do not induce much error. 

We would like to emphasize that the goal of the current paper is to propose a framework for designing a deterministic solver for the proton distribution. Fully parallelizing the code and comparing its efficiency with well-developed Monte Carlo (MC) codes are beyond the scope of the present work. However, there are several avenues for extending this work.
%Firstly, to achieve a sufficiently accurate proton dose distribution, catastrophic scatter interactions cannot be ignored. Although such scattering cannot be modeled with a Fokker-Planck type of approximation, the primary proton distribution can serve as the source for proton-nuclear scattering. By combining the contributions from both primary and secondary protons, one can obtain a much more accurate dose distribution. This is part of our ongoing work.

Since the depth variable $x$ is analogous to the time variable in classical time evolutionary problems, extending to a higher-order scheme is straightforward. This can be achieved by employing higher-order splitting methods, as demonstrated in \cite{cervi2020towards}. This will be the focus of our future work.

\bigskip
\textbf{Acknowledgement:} 
X. J. Zhang and M. Tang  were
 supported by NSFC12031013, Shanghai pilot innovation project 21JC1403500 and the Strategic Priority Research Program of Chinese Academy of Sciences Grant No.XDA25010401; X. J. Zhang was supported by  China Postdoctoral Science Foundation No.2022M722107.
\appendix
  \renewcommand{\appendixname}{Appendix~\Alph{section}}
  \section{Some physics data}\label{apen01}
In order to solve equation \eqref{f-e}, the expression of physics variables: stopping power, straggling coefficient and momentum transfer cross section should be known. In this section, we will give some  data about the above variables.

\subsection{Stopping power}\label{apen01}
\begin{figure}[ht]
	\centering
	\includegraphics[width=0.6\textwidth]{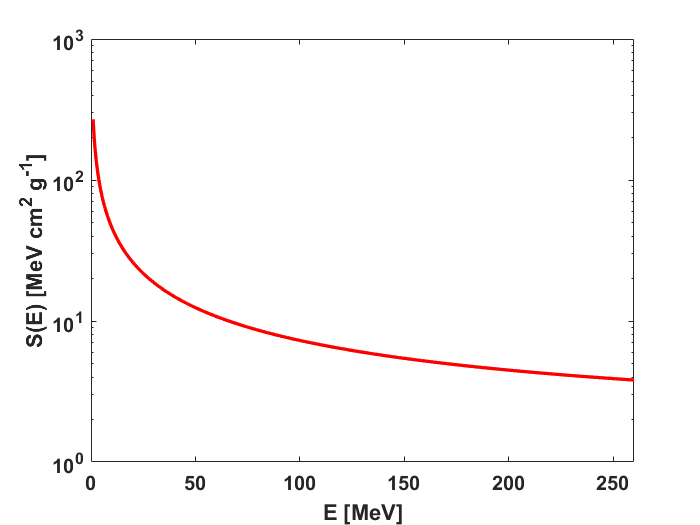}
		\caption{Stopping power of protons in water. }
		\label{fig:s_p}
	\end{figure}

The concept of energy loss has been explored by many researchers, and the  phenomena were described in different publications over the past years \cite{sigmund2004stopping,northcliffe1970range}. When the protons pass through the matter, the protons ionize atoms or molecules in their trajectory, leading to the energy loss of the protons. The stopping power describes the average loss of energy of the protons per unit of distance traveled, and  it models a property of the material, measured in $\text{MeV} \cdot \text{cm}^2 /\text{g}$. The stopping power is  described by the Bethe-Bloch formula \cite{groom2000passage}:
\begin{equation}
    S(E)=\frac{\kappa z_i^2Z}{\beta^2}\left(\text{ln}\frac{2m_ec^2\beta^2\gamma^2}{I}-\beta^2-\frac{\delta}{2}\right),
\end{equation}
where $\kappa$ is the stopping pre-factor (=0.307 $\text{MeV} \cdot \text{cm}^2 /\text{g}$), $z_i$ is the atomic number of incident particle, $Z$ is atomic number of medium material, and $\beta$ and $\gamma$ are the well knon relativistic quantities, which are defined as
$$
\beta=\frac{v}{c},\qquad \gamma=\sqrt{\frac{1}{1-\beta^2}},
$$
with $v$ the projectile's velocity and $c$ the speed of light. Moreover, $m_e$ is the electron mass, $\delta$ is the density correction term \cite{ziegler2010srim}, and $I$ is the ionisation potential \cite{turner1995atoms}, which is
\begin{align*}
\begin{split}
I= \left \{
\begin{array}{ll}
  19\quad \text{eV},                    & Z=1,\\
 11.2+11.7Z \quad \text{eV},     & 2\le Z\le 13,\\
 52.8+8.71Z\quad \text{eV},                                 & Z>13.
\end{array}
\right.
\end{split}
 \end{align*}
For a compound, its stopping power  can be represented as a convex combination of the individual stopping power of its constituent elements, i.e.,
$$
S(E)=\sum_k N_k S_k(E),
$$
in which the composite stopping power is obtained from the individual element $S_k(E)$ values weighted by $N_k$,  and the stopping power in water is plotted in Fig.~\ref{fig:s_p}.
\subsection{Straggling coefficient}\label{apen02}
When the protons undergo inelastic collisions with electrons, the overall direction of motion for the protons remains  unchanged due to the  mass of the protons is significantly greater than the electrons. However, energy loss occurs, accompanied by statistical fluctuations.This is because energy loss occurs through multiple collisions between the charged particles and the outer electrons of the target material, and the number of collisions as well as the energy transferred in each collision are randomly distributed. The average energy loss can be obtained by averaging over all particles, and the energy loss for individual particles is distributed around this average. This statistical distribution of energy loss is referred to as energy straggling \cite{PhysRevA.4.562}. 

\begin{figure}[ht]
	\centering
	\includegraphics[width=0.6\textwidth]{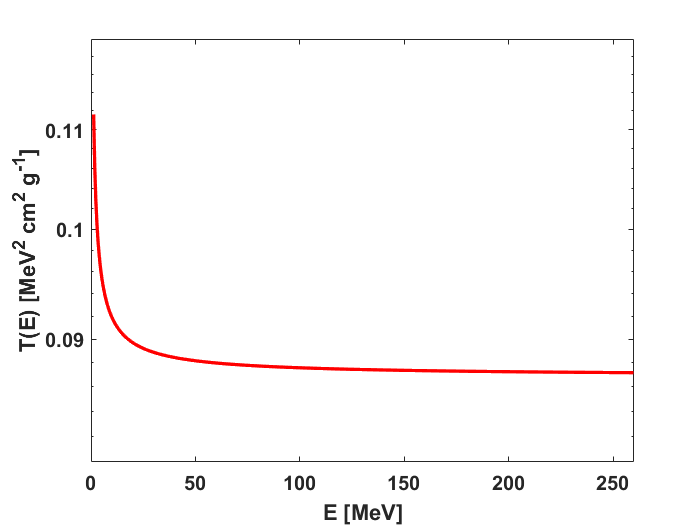}
		\caption{Straggling coefficient of protons in water. }
		\label{fig:s_c}
	\end{figure}

Regarding energy straggling, the earliest reference can be traced back to 1915. At that time, Bohr derived a formula for calculating the energy straggling of a single energetic charged particle passing through a medium based on the classical model of electrons \cite{sigmund2003binary}. In \cite{williams1932passage}, the author incorporated the quantum mechanical effects of electron binding in his derivation of the straggling coefficient, which is
\begin{equation}
    T(E)=4\pi e^4\left(1+\frac{4I}{3m_ev^2}\text{ln}\frac{2m_ev^2}{I}\right),
\end{equation}
with $e$ elementary charge, and the straggling coefficient in water is plotted in Fig.~\ref{fig:s_c}.

\subsection{Momentum transfer cross section}
\begin{figure}[ht]
	\centering
	\includegraphics[width=0.6\textwidth]{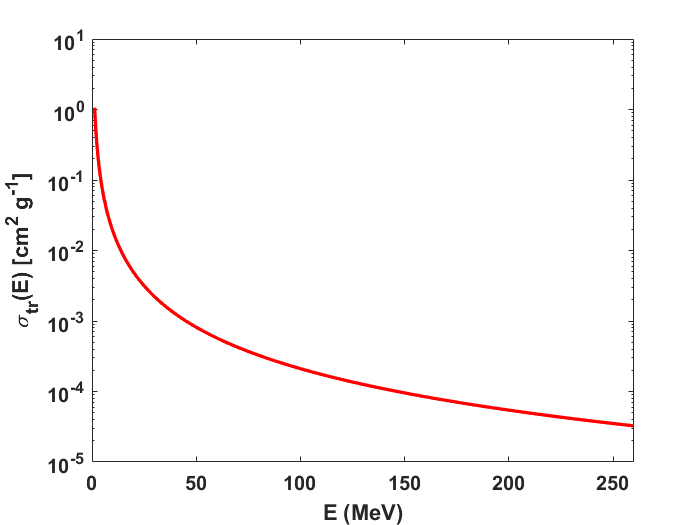}
		\caption{The momentum transfer cross section of protons in water. }
		\label{fig:m_t}
	\end{figure}

\begin{figure}[ht]
	\centering
	\includegraphics[width=0.6\textwidth]{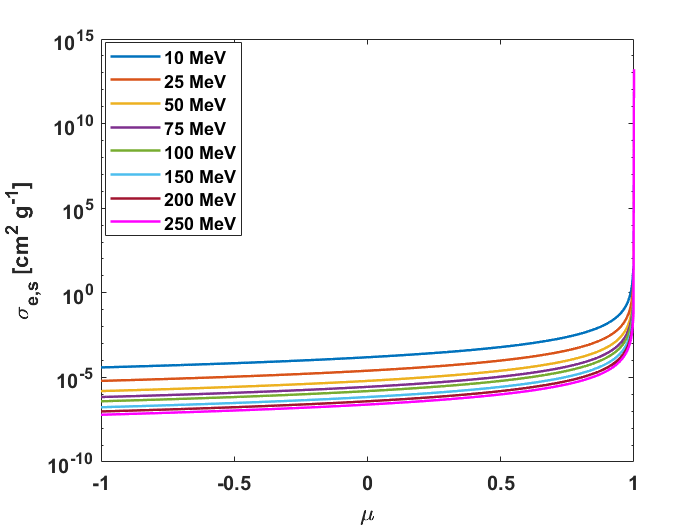}
		\caption{Screened Rutherford differential cross section of protons in water with different energies. }
		\label{fig:m_t1}
	\end{figure}

The incident protons undergo interactions with the nucleus through Coulomb forces, and these forces cause the incident protons to deviate from their initial trajectory. Therefore, in the elastic collisions between the protons and the nucleus,  the momentum transfer cross section is used to describe the angular diffusion. Recall the definition of the momentum transfer cross section:
\begin{equation*}
    \sigma_{tr}(E)=2\pi\int_{-1}^1(1-\mu)\sigma_{e,s}(\mu)d\mu.
\end{equation*}
Rutherford \cite{rutherford2012scattering} derived the differential cross sections for the elastic scattering process based on classical equations of motion, and  the Rutherford differential cross section 
$$
\sigma_{e,s} =\left( \frac{z_iZe^2}{4\pi\epsilon m_0v^2}\right)^2\frac{1}{(1-\mu)^2},
$$
is used to get the momentum transfer cross section, and $\epsilon$ is the vacuum permitivity, and the reduced mass of the proton is defined as
 $$
 \frac{1}{m_0}=\frac{1}{m}+\frac{1}{m_t},
 $$
 with the mass of proton $m$ and the ass of the target nucleus $m_t$. One can obtain the momentum transfer cross section:
 \begin{equation*}
    \sigma_{tr}(E)=2\pi\left( \frac{z_iZe^2}{4\pi\epsilon m_0v^2}\right)^2[\text{ln}(1-\mu)]|^{-1}_1.
\end{equation*}
However, the momentum transfer cross section which is gotten by using the Rutherford differential cross section becomes $\infty$ when $\mu\to 1$, so the result cannot be used directly. In practice, the screened Rutherford differential cross section \cite{goldstein2002classical,uilkema2012proton} is used:
% $$
% \sigma_{e,s} =\frac{\left(1+\frac{2\mu}{A_t}+\frac{1}{A_t^2}\right)^{\frac{3}{2}}}{1+\frac{\mu}{A_t}}\left( \frac{zZe^2}{4\pi\epsilon m_0v^2}\right)^2\frac{1}{(1-\mu+2\eta)^2},
% $$
 $$
\sigma_{e,s} =\left( \frac{z_iZe^2}{4\pi\epsilon m_0v^2}\right)^2\frac{1}{(1-\mu+2\eta)^2},
 $$
with 
$$
\eta = \Theta_{min}^2=\left(\frac{Z^{\frac{1}{3}}\alpha m_ec}{p}\right)^2,
$$
where $\alpha$ is the fine structure constant, $c$ is the speed of the light, and $p$ is the momentum of the proton, so the momentum transfer cross section is
 \begin{equation*}
    \sigma_{tr}(E)=2\pi\left( \frac{z_iZe^2}{4\pi\epsilon m_0v^2}\right)^2\left(\text{ln}\frac{\eta+1}{\eta}-\frac{1}{\eta+1}\right).
\end{equation*}
The momentum transfer cross section and the screened Rutherford differential cross section with different energies in water is plotted in Fig.~\ref{fig:m_t} and Fig.~\ref{fig:m_t1} respectively.
\subsection{Differential cross section for the catastrophic scattering}	\label{sec-d}
In this part, we will introduce a method to get the catastrophic scattering cross section based on data processing, and the data of the proton trajectory is from the MC software: FLUKA. Firstly, the differential cross section is constructed based on the discretization of the energy. And  for the energy discretization, the multi-group method is the most popular one \cite{mihalas1999foundations,pomraning2005equations}. In the multi-group method, the continuous energy space  is divided into $G$ groups, and each energy interval is denoted by $(E_{g-1/2},E_{g+1/2})$ ($g=1,\cdots, G$) with $g=G$ being the highest energy group. With the discretization of energy and the definition of catastrophic scattering differential and total cross section, in each energy group we have that
\begin{equation*}
    \sigma_{c,s,g}=\sigma_{c,t,g}\mathcal{P}\left(g'\to g,(u',v')\to(u,v)\right),
\end{equation*}
where $\mathcal{P}$ is the probability transition matrix. During the data of the secondary particles is less, the following assumption is adopted:
\begin{equation*}
    \mathcal{P}\left(g'\to g,(u',v')\to(u,v)\right)=\mathcal{P}^1\left(g'\to g\right)\mathcal{P}^2_g\left((u',v')\to(u,v)\right),
\end{equation*}
where $\mathcal{P}^1\left(g'\to g\right)$ is the probability that the protons in energy group $g'$ are scattered to the energy group $g$, and $\mathcal{P}^2_g\left((u',v')\to(u,v)\right)$ is the probability that the protons in direction $(u',v')$ are scattered to the direction $(u,v)$ in the energy group $g$. As in \cite{paganetti2002nuclear,fippel2004monte},  the dose from secondary protons shows a macroscopic build-up effect, so only secondary proton is considered in the catastrophic cross section.

The probability transition matrix is constructed by using the proton trajectory:
\begin{itemize}
\item Firstly, for energy probability transition matrix $\mathcal{P}^1\left(g'\to g\right)$, the protons can only lose energy in the catastrophic scatter interactions, so the protons move from higher energy to lower energy. In view of this reason, we use exponential distribution for fitting energy probability transition matrix $\mathcal{P}^1\left(g'\to g\right)$:
\begin{equation*}
\mathcal{P}^1\left(g'\to g\right)=\lambda e^{-\lambda(E_{g'}-E_g)}.
\end{equation*}
Here the fitting results in energy group 100, 200, 300 and 400 are displayed in Fig.~\ref{fig:2e}, where the energy domain [1 MeV, 260 MeV] is divided into 500 groups, and the injected energy is 230 MeV, and the target material is water.
\begin{figure}[h]
    \centering
   \subfloat{
		\includegraphics[width=0.5\textwidth]{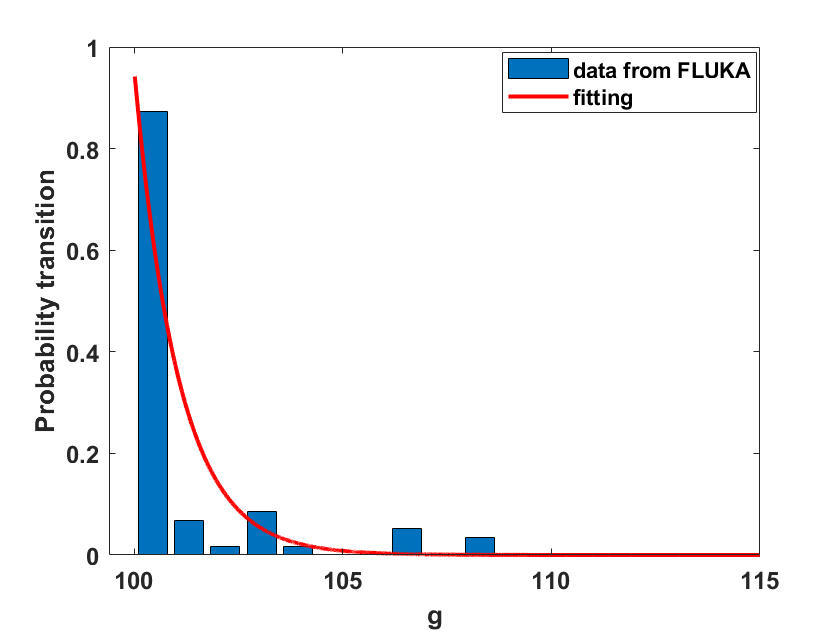}
	}
    \subfloat{
		\includegraphics[width=0.5\textwidth]{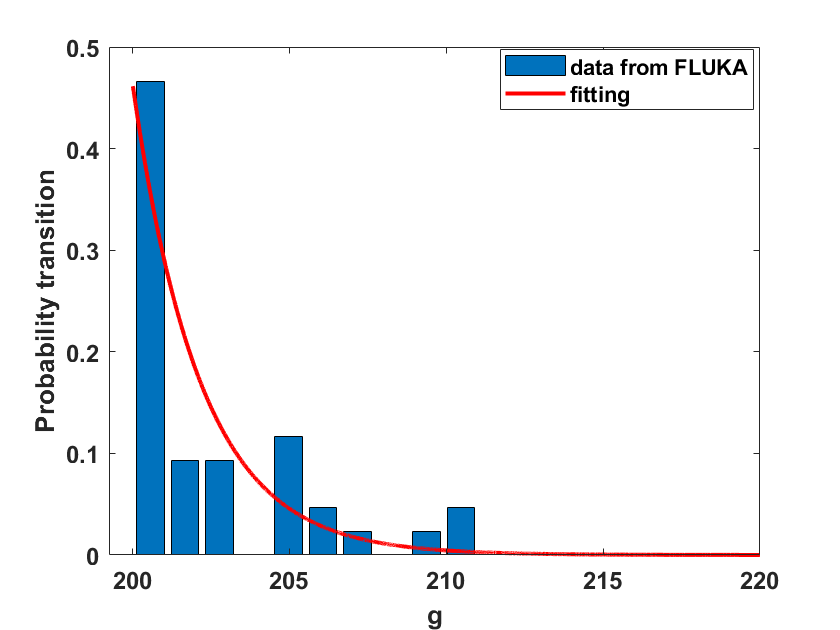}
	}
 \newline
    \subfloat{
		\includegraphics[width=0.5\textwidth]{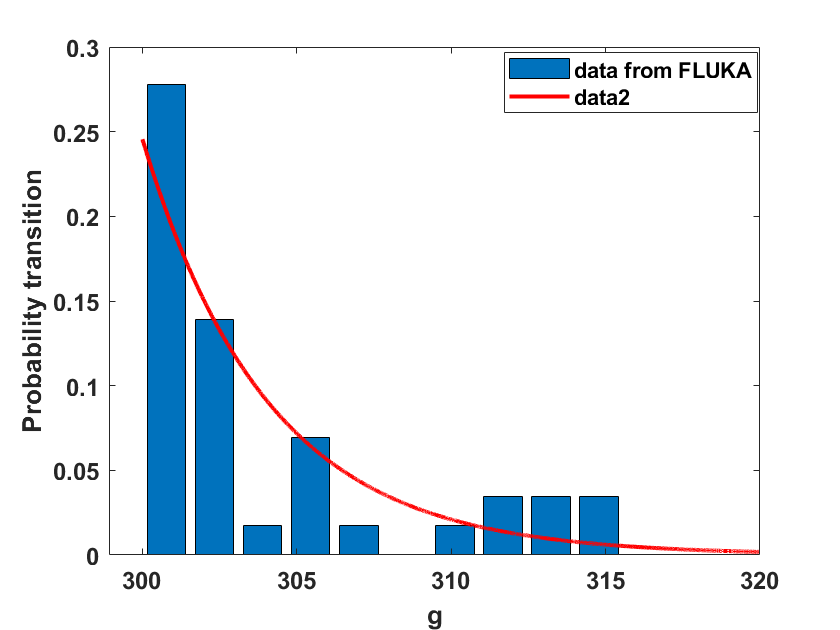}
	}
    \subfloat{
		\includegraphics[width=0.5\textwidth]{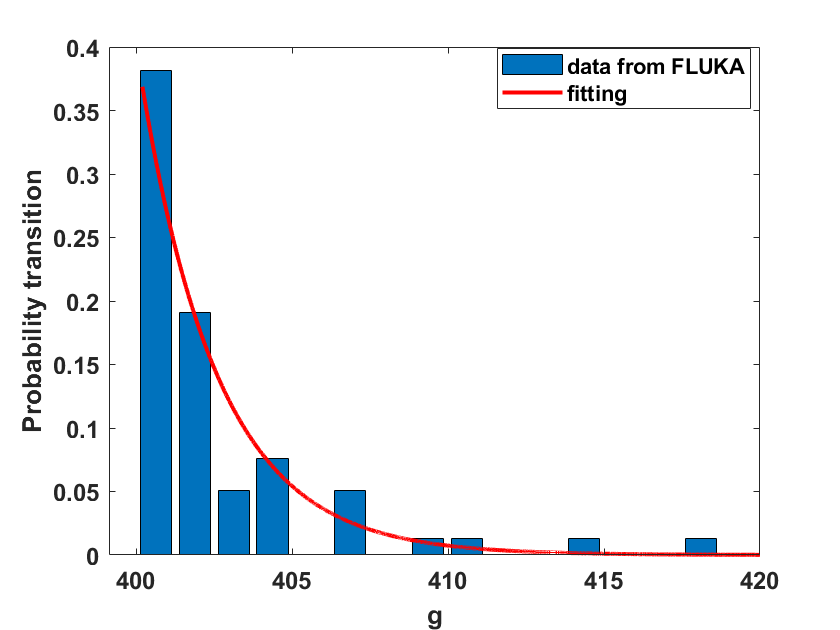}
	} 
    \caption{The energy probability transition fitting results in the energy group 100, 200, 300 and 400. }
    \label{fig:2e}
\end{figure}
% \begin{figure}[!h]
% 	\centering
% 	\subfigure[1]{
% 		\includegraphics[width=3.2in]{s_tr.png}
% 	}
% 	\subfigure[2]{
% 		\includegraphics[width=3.2in]{s_tr.png}
% 	}\\
% 	\subfigure[3]{
% 		\includegraphics[width=3.2in]{s_tr.png}
% 	}
% 	\subfigure[4]{
% 		\includegraphics[width=3.2in]{s_tr.png}
% 	}
% 	\caption{ evolve. }
% 	\label{fig.2_sgt}
% \end{figure}
\item For the fitting of the probability transition matrix regarding angles, since we can only obtain information about the deflected angles $\theta$ during collisions in FLUKA data, we need to follow several steps to obtain the probability transition matrix: 1. Fitting of the deflected angles; 2. Transformation of coordinates to the $(u, v)$ plane; 3. Normalization. 

Due to the fact that the cross section is forward-peaked, in this paper we use the Gamma distribution to fit the deflected angles during collisions:
\begin{equation*}
\mathcal{P}^2\left(\theta\right)=\frac{\beta^\alpha}{\Gamma
(\alpha)}\theta^{\alpha-1}e^{-\beta\theta},
\end{equation*}
Here the fitting results in energy group 100, 200, 300 and 400 are displayed in Fig.~\ref{fig:3a}.

We make the assumption that the medium through which the protons propagate is isotropic. In isotropic media, the scattering angle of the proton is independent of the angle at which the proton approaches the material. Under this assumption, the differential cross-section depends only on the scattering angle $\boldsymbol{\Omega}\cdot\boldsymbol{\Omega'}=\text{cos}\theta$, rather than on both $\boldsymbol{\Omega}$ and $\boldsymbol{\Omega'}$. By using the following equation and the above assumption,
\begin{equation*}
    \boldsymbol{\Omega}=(\mu,\eta,\xi),\quad u=\frac{\eta}{\mu},\quad v=\frac{\xi}{\mu},
    \end{equation*}
we have that
\begin{equation*}
\text{cos}\theta=\mu\mu'+\eta\eta'+\xi\xi'=\mu\mu'(1+uu'+vv'),
    \end{equation*}
moreover, by using $\mu^2+\eta^2+\xi^2=1$, we know that
\begin{equation*}
    \mu^2+\mu^2u^2+\mu^2v^2=1\Rightarrow\mu=(1+u^2+v^2)^{-\frac{1}{2}}.
\end{equation*}
From the above calculation, one can get
\begin{equation*}
    \theta=\arccos\left[\frac{1+uu'+vv'}{(1+u^2+v^2)^{\frac{1}{2}}(1+(u')^2+(v')^2)^{\frac{1}{2}}}\right]:=\arccos \Theta.
\end{equation*}

The final step is to normalize the obtained probability transition matrix to ensure the property of the sum of probabilities being 1. Then we can the following probability transition matrix:
\begin{equation*}
\mathcal{P}^2_g\left((u',v')\to(u,v)\right)=\frac{1}{K}\frac{\beta^\alpha}{\Gamma
(\alpha)}(\arccos \Theta)^{\alpha-1}e^{-\beta(\arccos \Theta)},
\end{equation*}
where
$$
K=\int_{-\infty}^{+\infty}\int_{-\infty}^{+\infty}\frac{\beta^\alpha}{\Gamma
(\alpha)}(\arccos \Theta)^{\alpha-1}e^{-\beta(\arccos \Theta)}\text{d}u'\text{d}v'.
$$
\begin{figure}[h]
    \centering
   \subfloat{
		\includegraphics[width=0.5\textwidth]{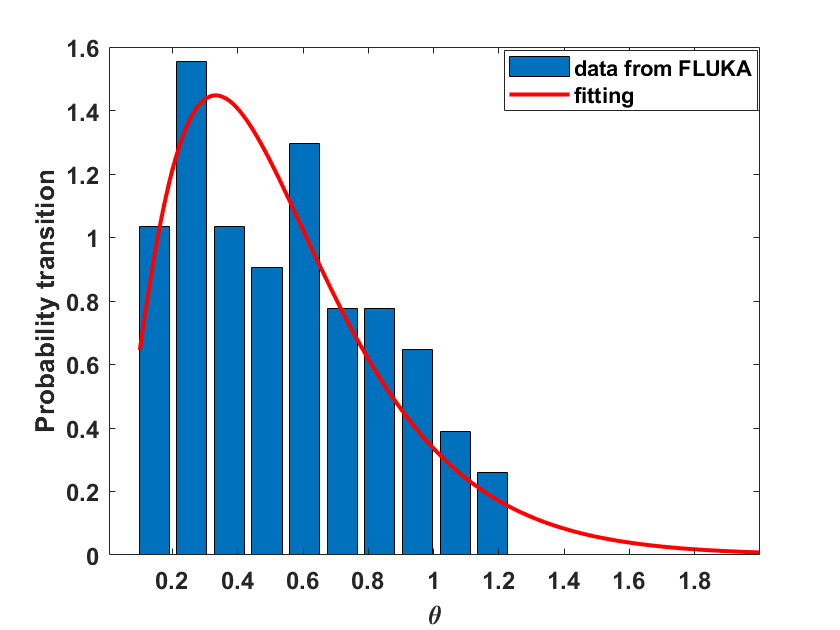}
	}
    \subfloat{
		\includegraphics[width=0.5\textwidth]{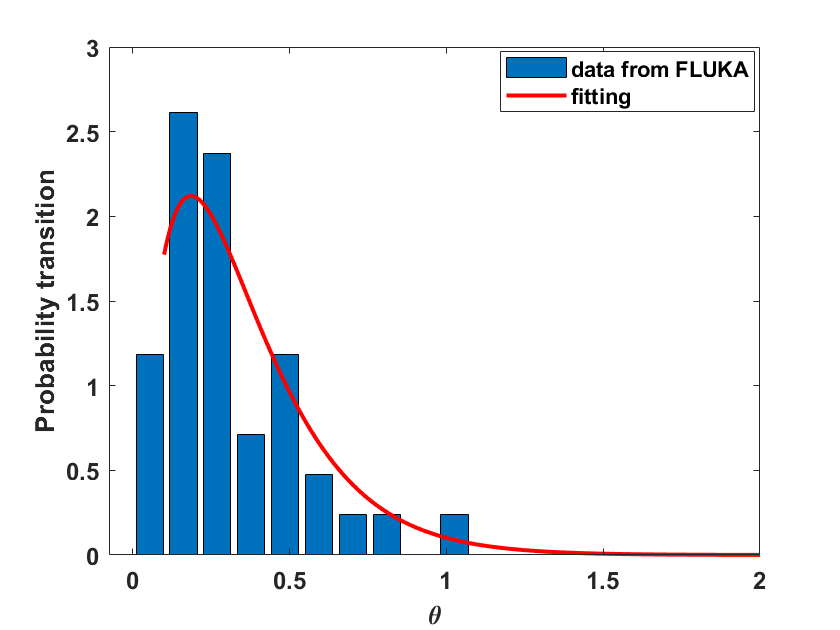}
	}
 \newline
    \subfloat{
		\includegraphics[width=0.5\textwidth]{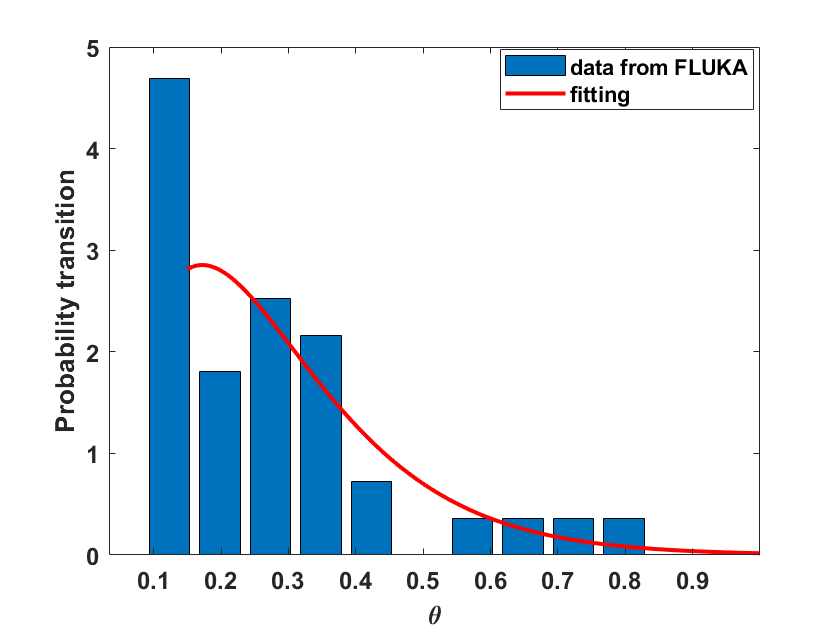}
	}
    \subfloat{
		\includegraphics[width=0.5\textwidth]{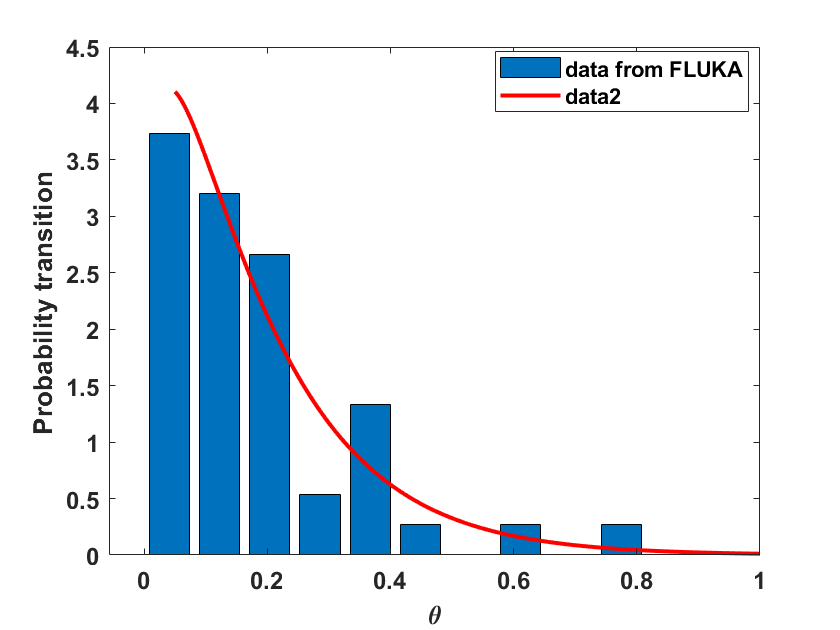}
	} 
    \caption{The angle probability transition fitting results in the energy group 100, 200, 300 and 400. }
    \label{fig:3a}
\end{figure}

\end{itemize}
\bibliography{sample}
\bibliographystyle{siam}

\end{document}